\documentclass[10pt]{amsart}
\usepackage{amsfonts, mathrsfs}
\usepackage{ifthen}
\usepackage{amsthm}
\usepackage{amsmath}
\usepackage{graphicx}
\usepackage{amscd,amssymb,amsthm}
\usepackage{graphicx}
\usepackage{epstopdf}
\usepackage{hyperref}
\usepackage{algorithm2e}
\usepackage{mathtools}
\usepackage{xcolor}

\newcounter{minutes}
\divide\time by 60
\newcounter{hours}
\multiply\time by 60 \addtocounter{minutes}{-\time}

\newtheorem{lemma}{Lemma}
\newtheorem{theorem}{Theorem}
\newtheorem{corollary}{Corollary}
\newtheorem{remark}{Remark}

\newcommand{\real}{\operatorname{Re}}

\newcommand{\logar}{\operatorname{Log}}

\numberwithin{equation}{section}


\begin{document}


\title[The radius of starlikeness of regular Coulomb wave functions]{The radius of starlikeness of regular Coulomb wave functions}

\author[\'A. Baricz]{\'Arp\'ad Baricz}
\address{Department of Economics, Babe\c{s}-Bolyai University, Cluj-Napoca 400591, Romania}
\address{Institute of Applied Mathematics, \'Obuda University, 1034 Budapest, Hungary}
\email{bariczocsi@yahoo.com}

\author[P. Kumar]{Pranav Kumar}
\address{Department of Mathematics,
	Indian Institute of Technology Madras, Chennai 600036, India}
\email{pranavarajchauhan@gmail.com}

\author[S. Singh]{Sanjeev Singh}
\address{Department of Mathematics,
	Indian Institute of Technology Indore, Indore 453552, India}
\email{snjvsngh@iiti.ac.in}

\def\thefootnote{}
\footnotetext{ \texttt{File:~\jobname .tex,
         printed: \number\year-\number\month-\number\day,
          \thehours.\ifnum\theminutes<10{0}\fi\theminutes}
} \makeatletter\def\thefootnote{\@arabic\c@footnote}\makeatother

\dedicatory{Dedicated to Mourad E.H. Ismail on the occasion of his 80th birthday}

\maketitle

\begin{abstract}
Motivated by the pioneering work of M.S. Robertson \cite{RobertSS} and R.K. Brown \cite{Brown,BrownUS}, who examined the geometric properties of some normalised solutions of second-order homogeneous differential equations, in this paper we investigate the radii of univalence and starlikeness for two kind of normalised regular Coulomb wave functions. Moreover, a generalized normalised Bessel function is introduced, and its radius of starlikeness is studied by using two different approaches. In addition, the asymptotic behaviour, with respect to the large order, of the radius of starlikeness of one type of normalised Coulomb wave functions is considered, which is in fact the first zero of the derivative of the regular Coulomb wave function. We derive a complete asymptotic expansion for this radius of starlikeness and provide a recurrence relation for the coefficients of this expansion. The proof is based on Rayleigh sums of the zeros of Coulomb wave functions, asymptotic inversion and some basic results on regular Coulomb wave functions developed by \v{S}tampach and \v{S}\'{t}ov\'{\i}\v{c}ek \cite{Flen_OrthPol_CWF}.
\end{abstract}

\section{\bf Introduction and Preliminaries}

In mathematics, a Coulomb wave function is a solution of the Coulomb wave equation, named after Charles-Augustin de Coulomb. The Coulomb wave functions are used to describe the behavior of charged particles in a Coulomb potential and can be written in terms of confluent hypergeometric functions or Whittaker functions of imaginary argument. The Coulomb differential equation reads
\begin{equation}\label{int_de}
W''(z)+\left(1-\frac{2\eta}{z}-\frac{l(l+1)}{z^2}\right)W(z)=0,
\end{equation}
where $z$ is the radial coordinate, $l$ is the angular momentum and $\eta$ is the Sommerfeld parameter. Two linearly independent solutions of \eqref{int_de} are $F_{l,\eta}(z) \text{ and } G_{l,\eta}(z)$, which are regular and irregular at $z=0$, respectively.  In accordance with \cite[eq.~33.2.4]{Nist} it is generally assumed that the standard regular Coulomb wave function $F_{l,\eta}(z) $ has a nonnegative real variable $z$, non-negative integer $l$, and a real parameter $\eta$. Nevertheless, various domains of theoretical physics, such as Regge pole theory and scattering theory necessitate the usage of complex variables $z$, $l$ and $\eta$ (see \cite{Humblet,Regge} and references therein). These are impetus behind the research endeavors of numerous scholars to investigate Coulomb wave functions utilizing complex variables and parameters (see \cite{Bar_Turen,BarFraHD,Alek_phaseintegral,Michel,Flen_OrthPol_CWF}). In order to contribute to this field the present study also undertakes an examination of the geometric properties of Coulomb wave functions involving complex variables. It is worth mentioning that the differential equation \eqref{int_de} is similar to that one studied by Nehari \cite{Nehari}, who considered the geometric properties of the solution of the differential equation of the form
\begin{equation}\label{int_de_1}
W''(z)+p(z)W(z)=0,
\end{equation}
where
\begin{equation}\label{p_exp}
z^2p(z)=p_0+p_1z+p_2z^2+\ldots
\end{equation}
is regular in $|z|<1$, whereby it has been demonstrated that, when $p_0=p_1=0$, no solution of \eqref{int_de_1} can assume a value zero more than once in $|z|<1$, given either
$$|p(z)|\leq\frac{1}{(1-|z|^2)^2}$$ or $$|p(z)|\leq\frac{\pi^2}{4}.$$
Later, Robertson \cite{RobertSS} generalized Nehari's result and investigated schlicht (univalent), starlike and spirallike properties of normalised solutions of \eqref{int_de_1} in the open unit disk. Subsequently, Brown \cite{Brown} extended Robertson's results for the solutions of \eqref{int_de_1} in the disk domain $|z|<r$ for $r>0$. Brown studied two normalised Bessel functions, specifically $[J_\nu(z)]^\frac{1}{\nu}$ and $z^{1-\nu}J_\nu(z)$, as examples of the above mentioned results. He was able to determine their radius of starlikeness and univalence for $\real \nu>0$ \cite[Theorem 2 and 3]{Brown}. Furthermore, Brown \cite{BrownUS} obtained the radius of starlikeness of one of the normalised Bessel function $z^{1-\nu}J_\nu(z)$, when $\nu\in (-\frac{1}{2}, 0)$. In \cite{RoSBesselBar}, the radius of starlikeness and convexity of normalised Bessel functions for $\nu>-1$ has been reconsidered. The authors in \cite{RoSBesselBar} used the Mittag-Leffler expansion and the distribution of zeros of Bessel functions of the first kind to perform this analysis. Note that, other geometric properties of Bessel functions of the first kind were discussed in details in \cite{KreyTod,Wilf,BarPunSC,BS_Besl_C,BarSanHBF}.
	
The regular Coulomb wave function is a generalization of the Bessel function of the first kind, and the relation between them is given by
\begin{equation}\label{CWF_Besl_rln}
F_{l-\frac{1}{2},0}(z)=\sqrt{\frac{\pi z}{2}}J_{l}(z),
\end{equation}
where $l, z \in \mathbb{C}$, such that $l\neq -1,-2,\ldots$, and $J_l$ stands for the Bessel function of the first kind of order $l$. A vast body of literature exists concerning the geometric properties and zeros of Bessel functions, whereas the research on the Coulomb wave function from the point of view of geometric function theory is comparatively scarce. Their close connection shows resemblance in their geometric properties and naturally motivates researchers to adopt similar approaches to study properties of Coulomb wave functions. However, it is interesting to note that the technique used in \cite{RoSBesselBar} for Bessel functions is not applicable for Coulomb wave function due to non-symmetric distribution of its zeros. All the same, motivated by the works of Hyden and Merkes \cite{HaydMerkChain} the authors \cite{BPS} worked in the direction of the starlikeness of regular Coulomb wave function where they find out the radii of disk, depending on parameters $l$ and $\eta$, which is mapped into a starlike domain by the normalised Coulomb wave functions. The authors in \cite{BPS} used the continued fraction expansion for the ratio ${zF'_{l,\eta}(z)}/{F_{l,\eta}(z)}$ and specialised those results also for Bessel functions. By using the idea of Brown \cite{Brown}, in this paper we find the exact radii of univalence and starlikeness of normalised Coulomb wave functions depending on some conditions on the parameters $l$ and $\eta$. This improves and complements the results obtained in \cite{BPS}. Note that this result is the best possible depending on the conditions of parameters and enhances what we obtained in \cite{BPS}.

An important study in this field is a recent paper of Baricz and Nemes \cite{BN21} in which the authors established the complete asymptotic expansions for the radii of starlikeness of two types of normalised Bessel functions of the first kind.  The Rayleigh sums associated with the zeros of Bessel functions of the firs kind, defined as sum of reciprocal of even powers of its zeros, play an important role to prove the main result of \cite{BN21}.  Due to the symmetrical distribution of the zeros of Bessel functions of the firs kind around origin, the Rayleigh sums associated with odd powers of its zeros vanish. In this paper, by using a similar approach as in \cite{BN21} we also determine the complete asymptotic expansion of the radius of starlikeness for one of the normalised Coulomb wave functions, however in the case of Coulomb wave functions the Rayleigh sums associated to odd powers do not vanish, which makes the approach more technical and comprehensive as compared to \cite{BN21}. It is also worth to mention here that some basic results on regular Coulomb wave functions developed by \v{S}tampach and \v{S}\'{t}ov\'{\i}\v{c}ek \cite{Flen_OrthPol_CWF} play an important role in the proofs of our main results related to regular Coulomb wave functions.

The structure of the paper is as follows: the subsequent portion of this section covers some basic definitions. Section 2 contains some well established results and the outcomes are related to the radii of univalence and starlikeness of normalised regular Coulomb wave functions. Furthermore, we provide some corollaries regarding Bessel functions. Section 3 is devoted to the study of the asymptotic behaviour, with respect to the large order, of the radius of starlikeness of one type of normalised Coulomb wave function. We derive a complete asymptotic expansion for the radius of starlikeness by using the asymptotic behaviour of the Rayleigh sums (convergent Laurent series expansions at infinity) and provide a recurrence relation for the coefficients of this expansion by using asymptotic inversion. Finally, in section 4, we provide the proofs of the preliminary and main results.

The conclusion of the paper is that indeed from the point of view of geometric function theory the Coulomb wave functions share similar properties to Bessel functions of the first kind and although the proofs in some cases are more technical in the case of Coulomb wave functions, the results resemble to the classical results on Bessel functions of the first kind. It is important to mention here that as a byproduct of our main results in this paper we obtain a complete asymptotic expansion for the smallest zero in modulus of the derivative of regular Coulomb wave function and this result may be useful in problems of applied mathematics or mathematical physics, where the zeros of the Coulomb wave function appear.

We would like to attract the attention to our previous work \cite{BPS}, where we investigated the starlikeness of Coulomb wave functions for positive Sommerfeld parameter $\eta$ by using continued fractions. In contrast, in this paper, we determine the radius of starlikeness of Coulomb wave functions for negative $\eta$ by using the technique of differential equations. Additionally, we present the asymptotic expansion for this radius of starlikeness when $\eta<0$. For the special case of $\eta=0$, the Coulomb wave function reduces to Riccati-Bessel function and the zeros of its derivatives were investigated by Boyer in \cite{Bo69}. Notably, the task of obtaining the radius of starlikeness of Coulomb wave functions for positive $\eta$ remains an open problem, with some progress outlined in \cite{BPS}. We anticipate that the results concerning the asymptotic expansion of the radius of starlikeness presented in this paper will also hold true for $\eta>0$ once the problem regarding radius of starlikeness of Coulomb wave functions is resolved for $\eta>0$.

Now, we recall some basic definitions, which will be used in the sequel. Specifically, for $r>0$, we denote the open disk of radius $r$ by $\mathbb{D}_r=\{z\in \mathbb{C}: |z|<r\}$. Moreover, we define the class of analytic functions $f:\mathbb{D}_r\to \mathbb{C}$, which satisfy the normalisation conditions $f(0)=0$ and $f'(0)=1$, i.e. they can be written in the form
\begin{equation}\label{normalized}
f(z)=z+\sum_{n\geq2}^{}a_nz^n,
\end{equation}
where the coefficients $a_n$ are real or complex numbers. We say that the function $f$, defined by \eqref{normalized}, is $\kappa\text{-spirallike}$, for some $|\kappa|<\frac{\pi}{2}$, in the disk $\mathbb{D}_r$ if $f$ is univalent in $\mathbb{D}_r$, and
$$\real\left(e^{i\kappa}\frac{zf'(z)}{f(z)}\right)>0\quad \mbox{for all}\quad z\in \mathbb{D}_r.$$
The function defined in \eqref{normalized} is called starlike in the disk $\mathbb{D}_r$ if $f$ is univalent in $\mathbb{D}_r$, and $f(\mathbb{D}_r)$ is a starlike domain in $\mathbb{C}$ with respect to the origin. These functions are characterized as
$$\real\left(\frac{zf'(z)}{f(z)}\right)>0\quad \mbox{for all}\quad z\in\mathbb{D}_r.$$
Note that starlike functions are special cases of spirallike functions when $\kappa=0$. For $\beta\in\left[0,1\right)$ we say that the function $f$ is starlike of order $\beta$ if and only if
$$\real\left(\frac{zf'(z)}{f(z)}\right)>\beta\quad \mbox{for all}\quad z\in\mathbb{D}_r.$$
The real number
$$r_\beta^*(f)=\sup\left\{r>0\left|\real\left(\frac{zf'(z)}{f(z)}\right.\right)>\beta\quad \mbox{for all}\quad z\in \mathbb{D}_r\right\}$$
is called the radius of starlikeness of order $\beta$ of the function $f$. We use the notation $r^*(f)=r_0^*(f)$ and this is the largest radius such that the image region $f(\mathbb{D}_{r^*(f)})$ is a starlike domain with respect to the origin. Similarly, the radius of univalence is the largest radius of the disk under which the function $f$ is univalent.

\section{\bf The radii of univalence and starlikeness of Coulomb wave functions}
\setcounter{equation}{0}
	
First we consider the regular Coulomb wave function, which is defined by (see \cite[eq.~33.2.4]{Nist})
$$F_{L,\eta}(z)=z^{L+1}e^{-iz}C_L(\eta) {}_1F_1(L+1-i\eta, 2L+2; 2iz)= C_L(\eta)\sum_{n\geq0}a_{L,n}z^{n+L+1},$$
where $L,\eta, z\in \mathbb{C}$, and $$C_L(\eta)=\frac{2^L e^{-\frac{\pi \eta}{2}}|\Gamma(L+1+i\eta)|}{\Gamma(2L+2)},$$
$$a_{L,0}=1,  \qquad a_{L,1}=\frac{\eta}{L+1}, \qquad a_{L,n}=\frac{2\eta a_{L,n-1}- a_{L,n-2}}{n(n+2L+1)},\qquad n\in \{2,3,\ldots\},$$
and ${}_1F_1$ stands for the Kummer confluent hypergeometric function. Two important normalised forms of regular Coulomb wave function are
\begin{equation}\label{Nor_CWF1}
f_{L,\eta}(z)= [C_L^{-1}(\eta)F_{L,\eta}(z)]^{\frac{1}{L+1}}
\end{equation}
and
\begin{equation}\label{Nor_CWF2}
g_{L,\eta}(z)= C_L^{-1}(\eta)z^{-L}F_{L,\eta}(z)=\sum_{n\geq 0}^{} a_{L,n}z^{n+1}.
\end{equation}
	
To investigate the Coulomb differential equation in this section our aim is to use the approach of Brown \cite{Brown}. For this we suppose that
\begin{equation}\label{p*}
z^2p^*(z)=\sum_{n\geq0}p^*_nz^n, \quad\quad p^*_0\leq \frac{1}{4},
\end{equation}
is regular for $|z|<r$ and real on the real axis. Given any non-negative constant $c$, we consider the second-order homogeneous differential equation
\begin{equation}\label{gen_diff_eqn}
W''(z)+\left[c\left(p^*(z)-\frac{p^*_0}{z^2}\right)+\frac{p^*_0}{z^2}\right]W(z)=0.
\end{equation}
Moreover, let $\alpha^*$ be the larger root of the characteristic equation associated with \eqref{gen_diff_eqn} and let
\begin{equation}\label{W_c}
W_c(z)=z^{\alpha^*}\sum_{n\geq0}a_n^*(c)z^n, \qquad a^*_0(c)=1,
\end{equation}
be the unique solution for $|z|<r$, corresponding to $\alpha^*$. The following theorem of Brown \cite{Brown} is the key tool in the proof of our main results of this section.
	
\begin{theorem}\label{theA}{\rm \cite[Theorem 1]{Brown}}
Let $z^2p(z)$, defined as in \eqref{p_exp}, be regular for $|z|<r$ and satisfy the inequality
\begin{equation}\label{p_cnd}
\real\left[e^{i\gamma} z^2p(z)\right]\leq \left[c[|z|^2p^*(|z|)-p^*_0]+p^*_0\right]\cos\gamma
\end{equation}
where $c\geq0,$ $|\gamma|\leq \pi/2$ and $z^2p^*(z)$ is defined in \eqref{p*}. With $p(z)$ chosen in this manner we define
\begin{align*}
W(z)=z^{\alpha}\sum_{n\geq0}a_n z^n, \qquad a_0=1,\quad |z|<r,
\end{align*}
to be the unique solution of \eqref{int_de_1} for $|z|<r$ corresponding to the root with larger real part of the associated characteristic equation. Let $W_c(z)$ be defined as in \eqref{W_c}. Then
\begin{equation}\label{gen_result}
\real\left(e^{i\gamma}\frac{zW'(z)}{W(z)}\right)\geq|z|\frac{W_c'(|z|)}{W_c(|z|)}\cos\gamma
\end{equation}
for all $|z|\leq \rho<r.$
\end{theorem}
Now we are going to state the main theorems of this section.
\begin{theorem}\label{Theorem2}
Let the complex number $L$ satisfy the inequalities $\real L>-1$ and  $|\arg(L+1))|<\frac{\pi}{4}$. Then for $\eta\leq0$, the normalised Coulomb wave function $z\mapsto f_{L,\eta}(z)$, as defined in \eqref{Nor_CWF1}, is regular, univalent and spirallike for $|z|<\widetilde{\rho}_{l,\eta,1},$ where $l(l+1)=\real\left[L(L+1)\right]$, $l>-1$ and $\widetilde{\rho}_{l,\eta,1}$ is the smallest root in modulus of the equation $r\mapsto F'_{l,\eta}(r)=0$. In particular, if for real $L>-1$ we have $l=L,$ then the radius of univalence is $\widetilde{\rho}_{L,\eta,1}$, and the radius of starlikeness of order $\beta\in\left[0,1\right)$ of the function $z\mapsto f_{L,\eta}(z)$ is the smallest root in modulus of the equation $rF_{L,\eta}'(r)-\beta(L+1)F_{L,\eta}(r)=0.$
\end{theorem}

\begin{remark}
{\em It is worth to mention that the theorem stated above remains valid for all values of $l$ that satisfy the provided conditions. Furthermore, the conditions $\real L>-1$ and $|arg(L+1)|<\frac{\pi}{4}$ imply $\real[L(L+1)]\geq-\frac{1}{4}$, thereby ensuring the existence of some real solution $l$. Similar observations hold for the subsequent theorems.}
\end{remark}

In particular, when $\beta=0,$ we have the following result.
\begin{corollary}
If $l>-1$ and $\eta\leq 0$, then the radius of starlikeness of $f_{l,\eta}(z)$ is $\widetilde{\rho}_{l,\eta,1}$, which is the smallest root in modulus of the equation $F'_{l,\eta}(r)=0$.
\end{corollary}
	
It is important to mention here that the relation \eqref{CWF_Besl_rln} suggests a result for Bessel functions corresponding to Theorem \ref{Theorem2}. For this consider the Legendre duplication formula for the Euler gamma function
$$\Gamma(2L+2)=\frac{1}{\sqrt{2\pi}}2^{2L+\frac{3}{2}}\Gamma(L+1)\Gamma\left(L+\frac{3}{2}\right).$$
By using the definition of $C_L(\eta)$ and the above equation we obtain that
\begin{equation}\label{C_exp}
C_L^{-1}(0)=\frac{\Gamma(2L+2)}{2^L\Gamma(L+1)}=\frac{1}{\sqrt{\pi}}2^{L+1}\Gamma\left(L+\frac{3}{2}\right).
\end{equation}
The normalised Bessel function related to $z\mapsto f_{L,\eta}(z)$ is obtained by taking $\eta=0$ as follows
\begin{align*}
f_{L-\frac{1}{2},0}(z)&=[C_{L-\frac{1}{2}}^{-1}(0)F_{L-\frac{1}{2},0}(z)]^{\frac{1}{L+\frac{1}{2}}}\\
&=\left[\frac{1}{\sqrt{\pi}}2^{L+\frac{1}{2}}\Gamma(L+1)\sqrt{\frac{\pi z}{2}}J_{L}(z)\right]^{\frac{1}{L+\frac{1}{2}}}\\
&=\left[2^L\Gamma(L+1)\sqrt{z}J_L(z)\right]^{\frac{1}{L+\frac{1}{2}}}.
\end{align*}

\begin{figure}[t]
\centering
\includegraphics[width=0.6\textwidth]{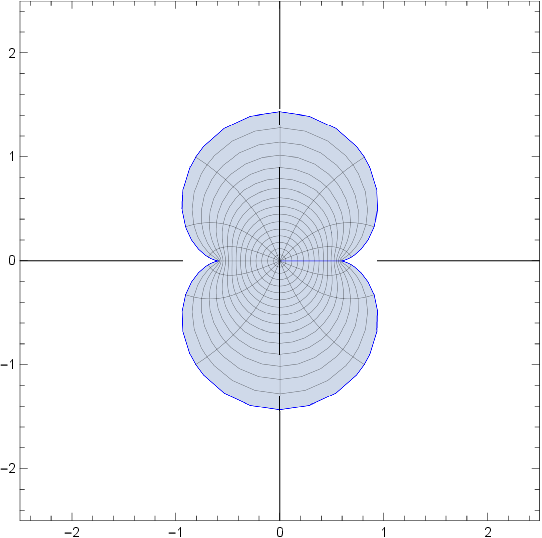}
\caption{The image of the open disk $\mathbb{D}_{\widetilde{\rho}_{0,0,1}}$ under the Bessel function $z\mapsto f_{-\frac{1}{2},0}(z)=\left[\sqrt{z}J_0(z)\right]^2,$ where $\widetilde{\rho}_{0,0,1}\sim 0.9407705639497375\ldots$ is the smallest positive zero of the function $r\mapsto 2rJ_0'(r)+J_0(r).$}
\label{Fig1}
\end{figure}

The next corollary is related to the function $z\mapsto f_{L-\frac{1}{2},0}(z)$. The case $L=0$ is illustrated in Figure \ref{Fig1}.

\begin{corollary}
Let the complex number $L$ satisfy the inequalities $\real L>-\frac{1}{2}$ and $\left|\arg\left(L+\frac{1}{2}\right)\right|<\frac{\pi}{4}$. Then the function $z\mapsto f_{L-\frac{1}{2},0}(z)=\left[2^L\Gamma(L+1)\sqrt{z}J_L(z)\right]^{\frac{1}{L+\frac{1}{2}}}$ is regular, univalent and spirallike for $|z|<\widetilde{\rho}_{l,0,1},$ where $l^2-\frac{1}{4}=\real\left(L^2-\frac{1}{4}\right)$, $l>-\frac{1}{2}$ and $\widetilde{\rho}_{l,0,1}$ is the smallest positive zero of the function $z\mapsto 2zJ'_l(z)+J_l(z).$ In particular, when $L>-\frac{1}{2}$, the radius of univalence is $\widetilde{\rho}_{L,0,1}$, and the radius of starlikeness of order $\beta\in\left[0,1\right)$ of the function $z\mapsto f_{L-\frac{1}{2},0}(z)$ is the smallest positive root of the equation $rJ'_{L}(r)-\left[\beta\left(L+\frac{1}{2}\right)-\frac{1}{2}\right]J_{L}(r)=0.$
\end{corollary}	

We mention that the normalised form $f_{L-\frac{1}{2},0}(z)$ is a novel addition to the existing studied normalised forms of Bessel functions. Motivated by the function $z\mapsto f_{L-\frac{1}{2},0}(z)$, we introduce a new generalized normalised Bessel function as
\begin{equation}\label{gen_nor_Bessel}
\varphi_{\nu,\alpha}(z)=\left[2^\nu\Gamma(\nu+1)z^\alpha J_\nu(z)\right]^{\frac{1}{\nu+\alpha}},
\end{equation}
for $\nu\neq-\alpha$ and $\alpha\in\mathbb{R}$. Note that
$$\varphi_{\nu,\alpha}(z)=\exp\left[\frac{1}{\nu+\alpha}\logar \left(2^\nu\Gamma(\nu+1)z^\alpha J_\nu(z) \right)\right],$$
where $\logar$ represents the principle branch of the logarithm function. Furthermore, all many-valued functions under consideration in this paper are evaluated with the principal branch.

The following theorem is related to the spirallike and starlike properties of the function $z\mapsto \varphi_{\nu,\alpha}(z)$. The proof of this theorem on the one hand is motivated by the results presented by Brown \cite{Brown}, and on the other hand an alternative proof for part {\bf b} is provided by using the idea of Baricz et al. \cite{RoSBesselBar}.

\begin{theorem}\label{the_Nor_Besl}
Let $z\mapsto \varphi_{\nu,\alpha}(z)$ denote the normalised Bessel function as defined in \eqref{gen_nor_Bessel}. Then the following assertions are true:
\begin{enumerate}
\item[\bf a.] If the complex number $\nu$ satisfies the inequalities $\real\nu>0,$ $|arg(\nu)|<\frac{\pi}{4}$, then $z\mapsto \varphi_{\nu,\alpha}(z)$ is regular, univalent and spirallike for $|z|<\rho_{\mu,\alpha},$ where $\mu^2=\real \nu^2,$ $\mu>0$ and $\rho_{\mu,\alpha}$ is the smallest positive root of the equation $\alpha J_\mu(r)+rJ_\mu'(r)=0$.
\item[\bf b.] If $\nu$  and $\nu+\alpha$ are positive, then the radius of starlikeness of order $\beta\in\left[0,1\right)$ of the function $z\mapsto \varphi_{\nu,\alpha}(z)$ is the smallest positive root of the equation $\left[\alpha-\beta(\nu+\alpha)\right]J_\nu(r)+rJ_\nu'(r)=0.$
\end{enumerate}
\end{theorem}

Note that Theorem \ref{the_Nor_Besl} is in fact the generalization of the result obtained in \cite[Theorem 1(a)]{RoSBesselBar}. It is worth also to mention that for $\alpha=0$, the above theorem yields results for the normalised Bessel function $\varphi_{\nu,0}(z)=\left[2^\nu\Gamma(\nu+1) J_\nu(z)\right]^{\frac{1}{\nu}}$, which has been extensively studied in the literature (see for example \cite{RoSBesselBar,BS_Besl_C} and the references therein). The introduction of this generalized normalised form $\varphi_{\nu,\alpha}(z)$ has new possibilities for improvements concerning the results related to Bessel functions. Specifically, when $\beta=0$, Theorem \ref{the_Nor_Besl} provides the radius of starlikeness for the function $z\mapsto \varphi_{\nu,\alpha}(z)$ as follows.

\begin{figure}[t]
\centering
\includegraphics[width=0.6\textwidth]{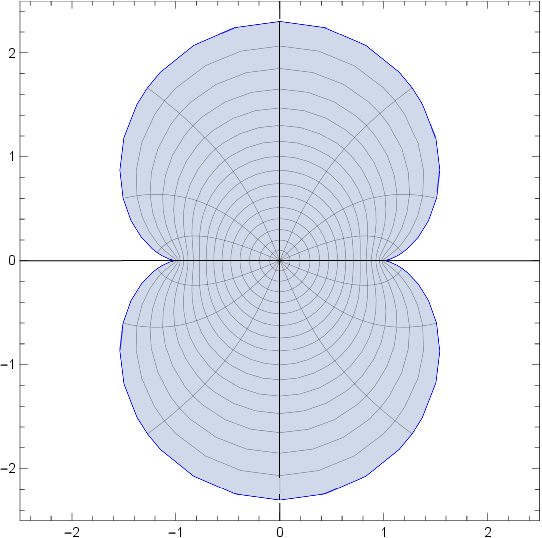}
\caption{The image of the open disk $\mathbb{D}_{\rho^*_{\frac{1}{2},0,1}}$ under the trigonometric function $z\mapsto g_{0,0}(z)=\sqrt{\frac{\pi z}{2}}J_{\frac{1}{2}}(z)=\sin z,$ where $\rho^*_{\frac{1}{2},0,1}\sim 1.5707963267948968\ldots$ is the smallest positive zero of the function $r\mapsto 2rJ_{\frac{1}{2}}'(r)+J_{\frac{1}{2}}(r).$}
\label{Fig2}
\end{figure}

\begin{corollary}
If $\nu>0$ and $\nu+\alpha>0$, then the radius of starlikeness of $z\mapsto \varphi_{\nu,\alpha}(z)$ is the smallest positive root of the equation $\alpha J_\nu(r)+rJ_\nu'(r)=0.$
\end{corollary}

Now, we focus on the normalized Coulomb wave function defined by \eqref{Nor_CWF2}.

\begin{theorem}\label{Theorem3}
Let the complex number $L=x+\mathrm{i}y$ satisfy the conditions:
\begin{equation}\label{condn}
x<1, \qquad y^2<x(x+1)+\frac{1}{4}.
\end{equation}
Then for $\eta\leq0$, the normalised Coulomb wave function $z\mapsto g_{L,\eta}(z)$, as defined in \eqref{Nor_CWF2}, is regular, univalent and spirallike for $|z|<\rho^*_{l,\eta,1},$ where $l(l+1)=\real\left[L(L+1)\right]$, $l>-1$ and $\rho^*_{l,\eta,1}$ is the smallest root in modulus of the equation $r\mapsto rF'_{l,\eta}(r)-\left(\real L\right)\cdot F_{l,\eta}(r)=0.$ In particular, when $L>-1$, then radius of univalence is $\rho^*_{l,\eta,1}$, and the radius of starlikeness of order $\beta\in\left[0,1\right)$ of the function $z\mapsto g_{L,\eta}(z)$ is the smallest root in modulus of the equation $rF'_{L,\eta}(r)-(\beta+L)F_{L,\eta}(r)=0.$
\end{theorem}

If $\beta=0$, then we arrive to the next result concerning the radius of starlikeness of $z\mapsto g_{l,\eta}(z).$

\begin{corollary}
If $l>-1$ and $\eta\leq 0$, then the radius of starlikeness of $z\mapsto g_{l,\eta}(z)$ is $\rho^*_{l,0,1}$, which is the smallest root in modulus of equation $rF'_{l,\eta}(r)-lF_{l,\eta}(r)=0.$
\end{corollary}

It is interesting to note that by using the equations \eqref{CWF_Besl_rln} and \eqref{C_exp}, we obtain the next normalised Bessel function related to $z\mapsto g_{L,\eta}(z)$ for $\eta=0$
$$g_{L-\frac{1}{2},0}(z)=2^L\Gamma(L+1)z^{1-L}J_L(z).$$
The radius of starlikeness of $z\mapsto g_{L-\frac{1}{2},0}(z)$ has been deduced in \cite[Theorem 1(b)]{RoSBesselBar} with the condition $L>-1$. Moreover, it is worth mentioning that the asymptotic behavior of the radius of starlikeness of the above function $z\mapsto g_{L-\frac{1}{2},0}(z)$ has been considered in \cite{BN21} in details. In addition, as a complementary contribution to the existing results, the next corollary concerns the univalence, spirallikeness and starlikeness of $z\mapsto g_{L-\frac{1}{2},0}(z)$. The case $L=\frac{1}{2}$ is illustrated in Figure \ref{Fig2}.

\begin{corollary}
Let the complex number $L=x+\mathrm{i}y$ satisfy the conditions $x<\frac{3}{2}$ and $y^2<x^2.$ Then $z\mapsto g_{L-\frac{1}{2},0}(z)=2^L\Gamma(L+1)z^{1-L}J_L(z)$ is regular, univalent and spirallike for $|z|<\rho^*_{l,0,1},$ where $l^2-\frac{1}{4}=\real\left[L^2-\frac{1}{4}\right],$ $l>-\frac{1}{2}$ and $\rho^*_{l,0,1}$ is the smallest positive zero of $r\mapsto rJ'_{l}(r)+\real(1-L)J_{l}(r)$. In particular, when $L$ is real and positive, then the radius of univalence is $\rho^*_{L,0,1}$ and the radius of starlikeness of order $\beta\in\left[0,1\right)$ of the function $z\mapsto g_{L-\frac{1}{2},0}(z)$ is the smallest positive root of the equation $rJ'_{L}(r)-(\beta+L-1)J_{L}(r)=0.$
\end{corollary}

\section{\bf Asymptotics of the radius of starlikeness of the normalised Coulomb wave function}
\setcounter{equation}{0}

In this section we consider the Rayleigh sums associated with the Coulomb wave function and its derivative as
\begin{equation}\label{Z_Coul}
	Z_\eta^{\left(k\right)}\left(L\right)=\sum_{n\geq 1}\frac{1}{\rho_{L,\eta,n}^k}
\end{equation}
and
\begin{equation}\label{Z_DCoul}
	\tilde{Z}_\eta^{\left(k\right)}\left(L\right)=\sum_{n\geq1}\frac{1}{\tilde{\rho}_{L,\eta,n}^k},
\end{equation}
where $\rho_{L,\eta,n}$ and $\tilde{\rho}_{L,\eta,n}$ are the non-zero roots of $z\mapsto F_{L,\eta}\left(z\right)$ and $z\mapsto \partial_z F_{L,\eta}\left(z\right)$, respectively. Here the zeros $\rho_{L,\eta,n}$ and $\tilde{\rho}_{L,\eta,n}$ are taken in such a way that $0<|\rho_{L,\eta,1}|<|\rho_{L,\eta,2}|<|\rho_{L,\eta,3}|<\ldots$ and $0<|\tilde{\rho}_{L,\eta,1}|<|\tilde{\rho}_{L,\eta,2}|<|\tilde{\rho}_{L,\eta,3}|<{\ldots}.$ Throughout this paper, if not stated otherwise, empty sums are taken to be zero. Additionally, $\mathbb{N}$ is the set of all positive integers and $\mathbb{N}_0:=\{0\}\cup \mathbb{N}$.

\begin{lemma}\label{lemma_asy1}
For $\alpha\in\mathbb{N}$ and $n\in \mathbb{N}_0$ consider the expression
\begin{equation}\label{p_n}
p_n^{\left(\alpha\right)}=\frac{(-1)^n}{2}\left(\frac{\alpha+1}{2}\right)^n.
\end{equation}
Then for any real $\eta,$ positive integer $k$ and real number $L>k+1$, the Rayleigh function $Z_\eta^{\left(2k\right)}\left(L\right)$ and $Z_\eta^{\left(2k+1\right)}\left(L\right)$ have the convergent Laurent series expansion
\begin{equation}\label{z_even_series}
Z_\eta^{\left(2k\right)}\left(L\right)=\frac{1}{L^{2k-1}}\sum_{n\geq0}\frac{\zeta_{n,\eta}^{\left(2k\right)}}{L^n}
\end{equation}
and
\begin{equation}\label{z_odd_series}
	Z_\eta^{\left(2k+1\right)}\left(L\right)=\frac{1}{L^{2k+1}}\sum_{n\geq0}\frac{\zeta_{n,\eta}^{\left(2k+1\right)}}{L^n}.
\end{equation}
Here the coefficients of the above Laurent series expansions are given by the following relations
$$\zeta_{0,\eta}^{\left(2\right)}=p_0^{\left(2\right)},\quad \zeta_{1,\eta}^{\left(2\right)}=p_1^{\left(2\right)}\quad \mbox{and}\quad \zeta_{n+2,\eta}^{\left(2\right)}=p_{n+2}^{\left(2\right)}+\sum_{m=0}^{n}\left(-1\right)^m\eta^2(m+1)p_{n-m}^{\left(2\right)}\quad \mbox{for}\quad n\in \mathbb{N}_0.$$
Moreover, for $k\in\{2,3,\ldots\}$ we have
$$\zeta_{0,\eta}^{\left(2k\right)}=\sum_{l=0}^{k-2}\zeta_{0,\eta}^{\left(2l+2\right)}\zeta_{0,\eta}^{\left(2k-2l-2\right)}p_0^{\left(2k\right)}, \quad \zeta_{1,\eta}^{\left(2k\right)}=\sum_{l=0}^{k-2}\sum_{q=0}^{1}\sum_{m=0}^{q}\zeta_{m,\eta}^{\left(2l+2\right)}\zeta_{q-m,\eta}^{\left(2k-2l-2\right)}p_{1-q}^{\left(2k\right)}$$
and
\begin{align*}\zeta_{n+2,\eta}^{\left(2k\right)}=&\sum_{l=1}^{k-2}\sum_{q=0}^{n}\sum_{m=0}^{q}\zeta_{m,\eta}^{\left(2l+1\right)}\zeta_{q-m,\eta}^{\left(2k-2l-1\right)}p_{n-q}^{2k} +\sum_{l=0}^{k-2}\sum_{q=0}^{n+2}\sum_{m=0}^{q}\zeta_{m,\eta}^{\left(2l+2\right)}\zeta_{q-m,\eta}^{\left(2k-2l-2\right)}p_{n+2-q}^{2k}\\
&+2\eta\sum_{l=0}^{n}\sum_{m=0}^{l}\left(-1\right)^mp_{l-m}^{\left(2k\right)}\zeta_{n-l,\eta}^{\left(2k-1\right)}
\end{align*}
for $n\in \mathbb{N}_0.$
Similarly, we have
$$\zeta_{n,\eta}^{\left(3\right)}=2\eta  \sum_{l=0}^{n}\sum_{m=0}^{l} \left(-1\right)^m p_{l-m}^{\left(3\right)}\zeta_{n-l,\eta}^{\left(2\right)}$$
and for $k\in\{2,3,\ldots\}$ the expression $\zeta_{n,\eta}^{\left(2k+1\right)}$ is given by the recurrence relation
\begin{align*}
\zeta_{n,\eta}^{\left(2k+1\right)}&=2\eta\sum_{l=0}^{n}\sum_{m=0}^{l}\left(-1\right)^m p_{l-m}^{\left(2k+1\right)}\zeta_{n-l,\eta}^{\left(2k\right)}+\sum_{l=1}^{k-1}\sum_{q=0}^{n}\sum_{m=0}^{q}\zeta_{m,\eta}^{\left(2l+1\right)}\zeta_{q-m,\eta}^{\left(2k-2l\right)}p_{n-q}^{\left(2k+1\right)}\\
&+\sum_{l=0}^{k-2}\sum_{q=0}^{n}\sum_{m=0}^{q}\zeta_{m,\eta}^{\left(2l+2\right)}\zeta_{q-m,\eta}^{\left(2k-2l-1\right)}p_{n-q}^{\left(2k+1\right)}
\end{align*}
for $n\in \mathbb{N}_0.$
\end{lemma}

Recall that in the case when $\eta\leq 0$ it is established in the previous section that the radius of starlikeness of the normalised Coulomb wave function is determined as the smallest root in modulus of equation involving the derivative of the function denoted as $F_{L,\eta}^\prime(z)$ or $\partial_z F_{L,\eta}\left(z\right).$ The next lemma will provide bounds for the radius of starlikeness of the Coulomb
wave function.

\begin{lemma}\label{lemma_asy2}
Let $L>-1$ and $\eta< 0$. For $s\in \mathbb{N}$ the radius of starlikeness $r^*\left(f_{L,\eta}\right)$ of the function $z\mapsto f_{L,\eta}\left(z\right)$ satisfies the inequalities
\begin{equation}\label{ros_ineqn}
\left(\tilde{Z}_\eta^{\left(2s\right)}\left(L\right)\right)^{-\frac{1}{s}}<\left(r^*\left(f_{L,\eta}\right)\right)^2<\frac{\tilde{Z}_\eta^{\left(2s\right)}\left(L\right)}{\tilde{Z}_\eta^ {\left(2s+2\right)}\left(L\right)},
\end{equation}
where
\begin{equation}\label{zeta_bar_exp}
\left\{\begin{array}{ll}
(2L+3)\tilde{Z}_\eta^{\left(2\right)}\left(L\right)&={1-La_1-\tilde{p}a_0+\tilde{p}^2},\\
(2L+4)\tilde{Z}_\eta^{\left(3\right)}\left(L\right)&={-La_2-\tilde{p}a_1+a_0\tilde{Z}_\eta^{\left(2\right)}\left(L\right)-2\tilde{p}\tilde{Z}_\eta^{\left(2\right)}\left(L\right)},\\
(2L+n+5)\tilde{Z}_\eta^{\left(n+4\right)}\left(L\right)&=-La_{n+3}-\tilde{p}a_{n+2}+\sum\limits_{m=0}^{n+1}a_m\tilde{Z}_\eta^{\left(3+n-m\right)}\left(L\right)\\
&\quad\quad\quad\quad+\sum\limits_{m=0}^{n}\tilde{Z}_\eta^{\left(m+2\right)}\left(L\right)\tilde{Z}_\eta^{\left(n-m+2\right)}
\left(L\right)-2\tilde{p}\tilde{Z}_\eta^{\left(n+3\right)}\left(L\right),
\end{array}\right.
\end{equation}
for $n\in \mathbb{N}_0$ and
$\tilde{p}=\frac{\left(L+2\right)}{\left(L+1\right)^2}\eta$. Moreover, the coefficients $a_n$ are given by the generating function
$$\frac{2\left(\rho-\eta\right)}{\rho^2-2\eta\rho-L\left(L+1\right)}=\sum_{n\geq0}a_n\rho^n,$$ where
\begin{equation}\label{coeff_a_n}
	a_0=\frac{2\eta}{L(L+1)},\quad a_1=-\frac{2(L^2+L+2\eta^2)}{L^2(L+1)^2}, \quad a_2=\frac{2\eta\left(3L^2+3L+4\eta^2\right)}{L^3(L+1)^3},\quad {\ldots}.
\end{equation}
In particular for $s=1$ we have
$$\left(\tilde{Z}_\eta^{\left(2\right)}\left(L\right)\right)^{-1}<\left(r^*\left(f_{L,\eta}\right)\right)^2<\frac{\tilde{Z}_\eta^{\left(2\right)}\left(L\right)}{\tilde{Z}_\eta^{\left(4\right)}\left(L\right)},$$
where
$$(2L+5)\tilde{Z}_\eta^{\left(4\right)}\left(L\right)={-La_3-\tilde{p}a_{2}+\sum_{m=0}^{1}a_m\tilde{Z}_\eta^{\left(3-m\right)}\left(L\right)+\tilde{Z}_\eta^{\left(2\right)}
\left(L\right)\tilde{Z}_\eta^{\left(2\right)}\left(L\right)-2\tilde{p}\tilde{Z}_\eta^{\left(3\right)}\left(L\right)}.$$
\end{lemma}

Observe that the previous lemma provides in fact some bounds for the square of the first zero $\tilde{\rho}_{L,\eta,1}$ of the derivative of the Coulomb wave functions, that is, for $L>-1,$ $\eta<0$ and $s\in \mathbb{N}$ the first zero $\tilde{\rho}_{L,\eta,1}$ satisfies the inequalities
\begin{align*}\label{}
	\left(\tilde{Z}_\eta^{\left(2s\right)}\left(L\right)\right)^{-\frac{1}{s}}<\left(\tilde{\rho}_{L,\eta,1}\right)^2<\frac{\tilde{Z}_\eta^{\left(2s\right)}\left(L\right)}{\tilde{Z}_\eta^ {\left(2s+2\right)}\left(L\right)},
\end{align*}
where for $k\in\{2,3,4,\ldots\}$ the expression $\tilde{Z}_\eta^{(k)}(L)$ is given by the equations in \eqref{zeta_bar_exp}. This result is an extension of the work of \v{S}tampach and \v{S}\'{t}ov\'{\i}\v{c}ek \cite{Flen_OrthPol_CWF} related to the Rayleigh sums associated with the derivative of the Coulomb wave function.

The next two results are important in the proof of the main result of this section, that is, Theorem \ref{asym_the}.

\begin{lemma}\label{Lemma3}
For $\eta\leq0,$ positive integer $k$ and large $L$, the Rayleigh sums $Z_\eta^{\left(2k\right)}\left(L\right)$ and $Z_\eta^{\left(2k+1\right)}\left(L\right)$ can be written as
\begin{equation}\label{even_zeta}
	Z_\eta^{\left(2k\right)}\left(L\right)=\frac{1}{L^{2k-1}}\mathcal{O}\left(1\right)
\end{equation}
and
\begin{equation}\label{odd_zeta}
Z_\eta^{\left(2k+1\right)}\left(L\right)=\frac{1}{L^{2k+1}}\mathcal{O}\left(1\right).
\end{equation}
\end{lemma}

\begin{lemma}\label{Lemma4}
For $\eta<0$ and large $L$ the radius of starlikeness $r^*\left(f_{L,\eta}\right)$ of the function $z\mapsto f_{L,\eta}\left(z\right)$
can be expressed as
$$r^*\left(f_{L,\eta}\right)=L\left(c+\mathcal{O}\left(\frac{1}{L}\right)\right),$$
where $c$ is some positive constant.
\end{lemma}

Before we state the next theorem, it is important to mention the concept of ordinary potential polynomials, as these polynomials assume a pivotal role in the proof of the subsequent theorem. Let $f(z)=1+\sum_{n\geq1}a_nz^n$ be a formal power series. Corresponding to $f(z)$, for any complex number $\alpha$, the ordinary potential polynomial $A_{\alpha,n}(\alpha_1,\alpha_2,\ldots,\alpha_n)$ is defined by generating function
$$\left(f(z)\right)^\alpha=\left(1+\sum_{n\geq0}\alpha_nz^n\right)^\alpha=\sum_{n\geq0}A_{\alpha,n}(\alpha_1,\alpha_2,\ldots,\alpha_n)z^n.$$
Thus, specifically $A_{\alpha,0}=1,$ $A_{\alpha,1}=\alpha \alpha_1$ and $A_{\alpha,2}=\alpha \alpha_2+\binom{\alpha}{2} \alpha_1^2$. One can refer to \cite[Appendix]{Nemes_LaplaceMethod} for additional details about the ordinary potential polynomials.

Finally, we state the main result of this section.

\begin{theorem}\label{asym_the}
For $k\in\{2,3,\ldots\}$ let $\zeta_{n,\eta}^{\left(k\right)}$ denote the coefficients as in Lemma \ref{lemma_asy1}. Then, for finite $\eta$ and $\eta<0$, the radius of starlikeness $r^*\left(f_{L,\eta}\right)$ has the following asymptotic expansion
\begin{align*}
	r^*\left(f_{L,\eta}\right) \sim L \left(c+\sum_{k\geq1} \frac{\epsilon_k}{L^k}\right)=L\left(\sqrt{2}+ \left(\sqrt{2}\eta +\frac{1}{2\sqrt{2}}-\frac{1}{2}\right)\frac{1}{L}+ \left(\frac{29}{16}+\frac{1}{4\sqrt{2}}-\frac{3}{4}\eta\right)\frac{1}{L^2}+\cdots \right)
\end{align*}	
as $L\to \infty$, where $c=\sqrt{2},$ $$2c\zeta_{0,\eta}^{\left(2\right)}\epsilon_1=c\eta-c^2\sum_{k=0}^{1}\left(-1\right)^{1-k}\zeta_{k,\eta}^{\left(2\right)}-\zeta_{0,\eta}^{\left(4\right)}A_{3,0}$$ and the coefficients $\epsilon_n$ for $n\in \mathbb{N}_0$ are given by the following recurrence relation
\begin{align}\label{rec_reln}
2c& \zeta_{0,\eta}^{(2)}\epsilon_{n+2}
=-2c\sum_{k=0}^{n}\epsilon_{k+1}\sum_{q=0}^{n-k+1}\left(-1\right)^{n-k-q+1}\zeta_{q,\eta}^{\left(2\right)}\\
&+\eta\sum_{k=0}^{n}\left(-1\right)^{n-k+1}\left(n-k\right)\epsilon_{k+1}+c^2\sum_{k=0}^{n+2}\left(-1\right)^{n-k+1}\zeta_{k,\eta}^{\left(2\right)}-
\sum_{j=0}^{n}\left(\sum_{k=0}^{n-j}\left(-1\right)^{n-j-k}\zeta_{k,\eta}^{\left(2\right)}\sum_{l=0}^{j}\epsilon_{l+1}\epsilon_{j-l+1}\right)\nonumber\\
&+\sum_{j=0}^{n+1}\left(-1\right)^{n-j}\sum_{m=2}^{j+2}\left(\sum_{k=0}^{j-m+2}\zeta_{j-m-k+2,\eta}^{\left(2m\right)}A_{m+1,k}\left(\epsilon_1,\ldots,\epsilon_k\right)\right)+\left(-1\right)^{n+1}c\eta\left(n+2\right)\nonumber\\
&+\sum_{j=0}^{n}\left(-1\right)^{n-j+1}\sum_{m=1}^{j+1}\left(\sum_{k=0}^{j-m+1}\zeta_{j-m-k+1,\eta}^{\left(2m+1\right)}A_{m+1,k}\left(\epsilon_1,\ldots,\epsilon_k\right)\right)\nonumber.
\end{align}
\end{theorem}

\section{\bf Proofs of the main results}
\setcounter{equation}{0}

\begin{proof}[\bf Proof of Theorem \ref{Theorem2}]
Let the complex number $L$ satisfy the inequalities $\real L>-1,$ $|\arg(L+1))|<\frac{\pi}{4}$ and suppose that $\eta\leq0$. For $c=1$ and
\begin{equation}\label{p*_l}
z^2p^*(z)=z^2-2\eta z-l(l+1),
\end{equation}
we write the differential equation \eqref{gen_diff_eqn} as
\begin{equation}\label{CWE_l}
W''(z)+\left(1-\frac{2\eta}{z}-\frac{l(l+1)}{z^2}\right)W(z)=0, \qquad |z|<r.
\end{equation}
Note that $-l(l+1)\leq \frac{1}{4}$ for all $l\in \mathbb{R}$, so \eqref{p*_l} is similar to \eqref{p*} with $p^*_0=-l(l+1),p^*_1=-2\eta,p^*_2=1$ and $p^*_i=0$ for $i\in\{3,4,5,\ldots\}$. Note that the condition $-l(l+1)\leq\frac{1}{4}$ complements $|\arg(L+1)|<\frac{\pi}{4}$. Further, let $W_1(z)$ be one of the solutions of \eqref{CWE_l} corresponding to the larger root of its associated characteristic equation $\lambda(\lambda-1)-l(l+1)=0$. We have
\begin{equation}\label{sol_w1}
W_1(z)=F_{l,\eta}(z),
\end{equation}
the regular Coulomb wave function with real parameters. Moreover, if we suppose that
\begin{align*}
z^2p(z)=z^2-2\eta z-L(L+1)
\end{align*}
then the equation \eqref{int_de_1} reduces to
\begin{equation}\label{CWE_2}
W''(z)+\left(1-\frac{2\eta}{z}-\frac{L(L+1)}{z^2}\right)W(z)=0, \qquad |z|<r.
\end{equation}
One of the solutions of \eqref{CWE_2}, corresponding to the root with larger real part, of the characteristic equation $\lambda(\lambda-1)-L(L+1)=0$ is
\begin{equation}\label{sol_w}
W(z)=F_{L,\eta}(z),
\end{equation}
the regular Coulomb wave function with complex parameter $L$. Note that the solutions \eqref{sol_w1} and \eqref{sol_w} are valid for all finite $z$. Furthermore, if $\eta\leq0$ and $l(l+1)=\real\left[L(L+1)\right],$ then
\begin{align*}
	\real\left[ z^2p(z)\right]&=\real z^2-2\eta \real z-\real\left[L(L+1)\right]\\
	&\leq |z|^2-2\eta|z|-l(l+1)\\
	&\leq |z|^2p^*(|z|).
\end{align*}
Therefore, for finite $z$ the condition \eqref{p_cnd} with $c=1$ and $\gamma=0$ is satisfied by $p(z)$ and $p^*(z)$. Consequently, in view of \eqref{gen_result} we have
\begin{equation}\label{special_result1}
\real\left(\frac{zW'(z)}{W(z)}\right)\geq|z|\frac{W'_1(|z|)}{W_1(|z|)}
\end{equation}
for all $|z|<r$. In addition, equations \eqref{sol_w1}, \eqref{sol_w} and \eqref{special_result1} imply that
\begin{equation}\label{final_reln}
\real\left(\frac{zF'_{L,\eta}(z)}{F_{L,\eta}(z)}\right)\geq|z|\frac{F'_{l,\eta}(|z|)}{F_{l,\eta}(|z|)}.
\end{equation}
Now, by taking the logarithmic derivative of both sides of \eqref{Nor_CWF1}, we arrive at
\begin{equation}\label{LogDer}
\frac{f_{L,\eta}'(z)}{f_{L,\eta}(z)}=\frac{1}{L+1}\frac{F_{L,\eta}'(z)}{F_{L,\eta}(z)}.
\end{equation}
If we compare the real parts after multiplying both sides of \eqref{LogDer} by $ze^{i\theta}$ and writing $L+1=re^{i\phi}$, for $|\phi|< \frac{\pi}{4}$, then we obtain
\begin{equation}\label{spiral_reln}
\real\left(e^{i\theta}\frac{zf'_{L,\eta}(z)}{f_{L,\eta}(z)}\right)=\frac{1}{r}\real\left(e^{i(\theta-\phi)}\frac{zF'_{L,\eta}(z)}{F_{L,\eta}(z)}\right).
\end{equation}
The minimum principle for harmonic functions together with \eqref{final_reln} and \eqref{spiral_reln}, shows that $z\mapsto f_{L,\eta}(z)$ is spirallike (for $\theta=\phi$) in disk $|z|=r<\widetilde{\rho}_{l,\eta,1}$, if $L$ is complex. Also, provided that $L$ is real, $r\mapsto F'_{L,\eta}(r)$ vanishes for $\widetilde{\rho}_{L,\eta,1}$, hence the function $z\mapsto f_{L,\eta}(z)$ for $L+1>0,$ cannot be univalent in any domain like $|z|=r>\widetilde{\rho}_{L,\eta,1}$.

If $L+1$ is real and positive then from \eqref{final_reln} and \eqref{LogDer} we obtain
$$\real\left(\frac{zf_{L,\eta}'(z)}{f_{L,\eta}(z)}\right)=\frac{1}{L+1}\real\left(\frac{zF_{L,\eta}'(z)}{F_{L,\eta}(z)}\right)\geq\frac{1}{L+1}\frac{rF_{L,\eta}'(r)}{F_{L,\eta}(r)}.$$
The above inequality along with the fact that $\lim_{z\to0} \left.{zf_{L,\eta}'(z)}\right/{f_{L,\eta}(z)}=1$ and the minimum principle for harmonic functions imply that
$$\real\left(\frac{zf_{L,\eta}'(z)}{f_{L,\eta}(z)}\right)>\beta$$
holds if and only if $|z|<\tau_{l,\eta,1}$ where $\tau_{l,\eta,1}$ is the smallest positive root of the equation
$$\frac{1}{L+1}\frac{rF_{L,\eta}'(r)}{F_{L,\eta}(r)}=\beta,$$
which is equivalent to $$rF_{L,\eta}'(r)-\beta(L+1)F_{L,\eta}(r)=0.$$
\end{proof}

\begin{proof}[\bf Proof of Theorem \ref{the_Nor_Besl}]
Let the complex number $\nu$ satisfy the inequalities  $\real\nu>0$ and $\left|\arg\nu\right|<\frac{\pi}{4}$. For $c=1$, $\mu^2=\real\nu^2$ and
\begin{equation}\label{p*_l_bsl}
z^2p^*(z)=z^2+\frac{1}{4}-\mu^2,
\end{equation}
we write the differential equation \eqref{gen_diff_eqn} as
\begin{equation}\label{NBF_1}
W''(z)+\left(1-\frac{1}{z^2}\left(\mu^2-\frac{1}{4}\right)\right)W(z)=0, \qquad |z|<r.
\end{equation}
Note that $\frac{1}{4}-\mu^2\leq \frac{1}{4}$, thus in this case \eqref{p*_l_bsl} is similar to \eqref{p*} with $p^*_0=\frac{1}{4}-\mu^2,$ $p_1^*=0,$ $p^*_2=1$ and $p^*_i=0$ for $i\in\{3,4,5,\cdots\}$. Let $W_1(z)$ be one of the solutions of \eqref{NBF_1} corresponding to the larger root of its associated characteristic equation $\lambda(\lambda-1)+\frac{1}{4}-\mu^2=0$. We have
\begin{equation}\label{sol_w1_bsl}
W_1(z)=2^\mu\Gamma(\mu+1)z^{1/2}J_\mu(z).
\end{equation}
Furthermore, suppose that
\begin{align*}
z^2p(z)=z^2+\frac{1}{4}-\nu^2,
\end{align*}
and then equation \eqref{int_de_1} reduces to
\begin{equation}\label{NBF_2}
W''(z)+\left(1-\frac{1}{z^2}\left(\nu^2-\frac{1}{4}\right)\right)W(z)=0, \qquad |z|<r.
\end{equation}
One of the solutions of \eqref{NBF_2} corresponding to roots, with larger real part, of the characteristic equation $\lambda(\lambda-1)+\frac{1}{4}-\nu^2=0$ is
\begin{equation}\label{sol_w_b}
W(z)=2^\nu\Gamma(\nu+1)z^{1/2}J_\nu(z),
\end{equation}
where $\nu$ is complex. Note that the solutions \eqref{sol_w1_bsl} and \eqref{sol_w_b} are valid for all finite $z$. Moreover, by using the conditions $\mu^2=\real\nu^2,$ $\mu>0$ we have
$$\real\left[z^2p(z)\right]=\real z^2+\frac{1}{4}-\real\nu^2\leq |z|^2+\frac{1}{4}-\mu^2=|z|^2p^*(|z|).$$
Therefore, for finite $z,$ $c=1$ and $\gamma=0,$ \eqref{p_cnd} is satisfied by $p(z)$ and $p^*(z)$. Consequently, \eqref{gen_result} implies
\begin{equation}\label{special_result}
\real\left(\frac{zW'(z)}{W(z)}\right)\geq|z|\frac{W'_1(|z|)}{W_1(|z|)}
\end{equation}
for all $|z|<r.$ Taking the logarithmic derivative of the corresponding parts of the equations \eqref{sol_w1_bsl} and \eqref{sol_w_b}, we obtain that
\begin{equation}\label{logd_wp}
\frac{W'_1(z)}{W_1(z)}=\frac{1}{2z}+\frac{J'_\mu(z)}{J_\mu(z)},
\end{equation}
and
\begin{equation}\label{logd_w}
\frac{W'(z)}{W(z)}=\frac{1}{2z}+\frac{J'_\nu(z)}{J_\nu(z)},
\end{equation}
respectively. In view of \eqref{special_result}, \eqref{logd_wp} and \eqref{logd_w}, we arrive at
\begin{equation}\label{Bessel_reln}
\real\left(\frac{zJ'_\nu(z)}{J_\nu(z)}\right)\geq\frac{|z|J'_\mu(|z|)}{J_\mu(|z|)}
\end{equation}
for all $|z|<r.$ Now, taking the logarithmic derivative of both sides of \eqref{gen_nor_Bessel}, we find that
\begin{equation}\label{LogDer_Bessel}
\frac{z\varphi_{\nu,\alpha}'(z)}{\varphi_{\nu,\alpha}(z)}=\frac{\alpha}{\nu+\alpha}+\frac{zJ_\nu'(z)}{(\nu+\alpha)J_\nu(z)}.
\end{equation}
By comparing the real parts after multiplying both sides of \eqref{LogDer_Bessel} by $ze^{\mathrm{i}\theta}$ and writing $\nu+\alpha=re^{\mathrm{i}\phi}$, for $|\phi|< \frac{\pi}{2}$, we obtain
\begin{equation}\label{spiral_reln_bessel}
\real\left(e^{\mathrm{i}\theta}\frac{z\varphi_{\nu,\alpha}'(z)}{\varphi_{\nu,\alpha}(z)}\right)=\real\left(\alpha e^{\mathrm{i}(\theta-\phi)} +\frac{e^{\mathrm{i}(\theta-\phi)}zJ_\nu'(z)}{J_\nu(z)}\right).
\end{equation}
The minimum principle for harmonic functions along with \eqref{Bessel_reln} and \eqref{spiral_reln_bessel}, suggests that if $\alpha+\mu>0$, the function $z\mapsto \varphi_{\nu,\alpha}(z)$ is spirallike (for $\theta=\phi$) for $|z|<\rho_{\mu,\alpha}$ where  $\rho_{\mu,\alpha}$ is the smallest positive zero of function $r\mapsto \alpha J_\mu(r)+rJ'_\mu(r)$. This completes the proof of part {\bf a.}

Now, let us assume that $\nu$ and $\nu+\alpha$ are real and positive. In this case equation \eqref{LogDer_Bessel} implies
\begin{equation}\label{Bes_alpha_real}
\real\left(\frac{z\varphi_{\nu,\alpha}'(z)}{\varphi_{\nu,\alpha}(z)}\right)=\frac{\alpha}{\nu+\alpha}+\frac{1}{\nu+\alpha}\real\left(\frac{zJ_\nu'(z)}{J_\nu(z)}\right).
\end{equation}
In view of \eqref{Bessel_reln}, \eqref{LogDer_Bessel} and \eqref{Bes_alpha_real} we arrive at
\begin{align*}
\real\left(\frac{z\varphi_{\nu,\alpha}'(z)}{\varphi_{\nu,\alpha}(z)}\right)&=\frac{\alpha}{\nu+\alpha}+\frac{1}{\nu+\alpha}\real\left(\frac{zJ_\nu'(z)}{J_\nu(z)}\right)\\
	&\geq\frac{\alpha}{\nu+\alpha}+\frac{|z|J_\nu'(|z|)}{(\nu+\alpha)J_\nu(|z|)}\\
	&=\frac{|z|\varphi_{\nu,\alpha}'(|z|)}{\varphi_{\nu,\alpha}(|z|)}.
\end{align*}
The above inequality along with the fact that $\lim_{z\to0} \left.{z\varphi_{\nu,\alpha}'(z)}\right/{\varphi_{\nu,\alpha}(z)}=1$ and the minimum principle for harmonic functions imply that
$$\real\left(\frac{z\varphi_{\nu,\alpha}'(z)}{\varphi_{\nu,\alpha}(z)}\right)>\beta$$
holds if and only if $|z|<\tau_{l,\alpha,\eta,1}$ where $\tau_{l,\alpha,\eta,1}$ is the smallest positive root of the equation
$$\frac{r\varphi_{\nu,\alpha}'(r)}{\varphi_{\nu,\alpha}(r)}=\beta,$$
which is equivalent to $$\left[\alpha-\beta(\nu+\alpha)\right]J_\nu(r)+rJ_\nu^\prime(r)=0.$$
\end{proof}

\begin{proof}[\bf Another proof of Theorem \ref{the_Nor_Besl}(b)]
Let $j_{\nu,n}$ denote the $n$-th positive zero of the Bessel function $z\mapsto J_\nu(z)$. The Weierstrassian decomposition of Bessel functions of the first kind is given by \cite[eq. 10.21.15]{Nist}
\begin{equation}\label{wess_decom}
J_\nu(z)=\frac{z^\nu}{2^\nu\Gamma(\nu+1)}\prod_{n\geq1}^{}\left(1-\frac{z^2}{j_{\nu,n}^2}\right),
\end{equation}
where the infinite product is uniformly convergent on each compact subset of $\mathbb{C}$. Logarithmic differentiation of both sides of \eqref{wess_decom} yields
\begin{equation}\label{Bes_logder_wess}
\frac{zJ_\nu'(z)}{J_\nu(z)}=\nu-\sum_{n\geq1}^{}\frac{2z^2}{j_{\nu,n}^2-z^2}.
\end{equation}
Moreover, from \eqref{LogDer_Bessel} and \eqref{Bes_logder_wess} we obtain that
$$\frac{z\varphi_{\nu,\alpha}'(z)}{\varphi_{\nu,\alpha}(z)}=\frac{\alpha}{\nu+\alpha}+\frac{zJ_\nu'(z)}{(\nu+\alpha)J_\nu(z)}=1-\frac{1}{\nu+\alpha}\sum_{n\geq1}^{}\frac{2z^2}{j_{\nu,n}^2-z^2}.$$
Let $x_{\nu,\alpha,\beta,1}$ denote the smallest positive root of the transcendental equation
$$\left[\alpha-\beta(\nu+\alpha)\right]J_\nu(r)+rJ_\nu'(r)=0.$$
We know that under the conditions $\alpha+\nu>0$ and $\nu>-1$, the Dini function, denoted as $z\mapsto zJ_\nu^\prime(z)+\alpha J_\nu(z)$ has only real zeros \cite[p. 597]{watson}. Moreover, based on the results of Ismail and Muldoon \cite[p. 11]{Ism_Muld}, we know that the smallest positive zero of the aforementioned function is smaller than $j_{\nu,1}$. Consequently, this indicates that $x_{\nu,\alpha,\beta,1}<j_{\nu,1}$ holds true for any $\nu>0$. Moreover, for all $0<\beta<1$ and $n\in\{2,3,4,\ldots\}$ we can deduce that $\mathbb{D}_{x_{\nu,\alpha,\beta,1}}$ $\subset$ $\mathbb{D}_{j_{\nu,1}}$ $\subset$ $\mathbb{D}_{j_{\nu,n}}$ when $\nu>0$. It is also known (see \cite{szasz} or \cite[eq. (2.4)]{RoSBesselBar}) that if $z\in\mathbb{C}$ and $\delta\in \mathbb{R}$ such that $\delta>|z|$, then
\begin{equation}\label{ineq}
\frac{|z|}{\delta-|z|}\geq\real\left(\frac{z}{\delta-z}\right).
\end{equation}
In view of \eqref{ineq}, for $\nu>-1,$ $n\in \{1,2,...\}$ and $z\in\mathbb{D}_{j_{\nu,1}}$ we have
\begin{equation}\label{ineq1}
\frac{|z|^2}{j^2_{\nu,n}-|z|^2}\geq\real\left(\frac{z^2}{j_{\nu,n}^2-z^2}\right).
\end{equation}
In view of the condition $\nu+\alpha>0$ and \eqref{ineq1} we conclude that
$$\real\left(\frac{z\varphi_{\nu,\alpha}'(z)}{\varphi_{\nu,\alpha}(z)}\right)=1-\frac{1}{\nu+\alpha}\real\left(\sum_{n\geq1}^{}\frac{2z^2}{j_{\nu,n}^2-z^2}\right)
\geq1-\frac{1}{\nu+\alpha}\frac{|z|^2}{j^2_{\nu,n}-|z|^2}=\frac{|z|\varphi_{\nu,\alpha}'(|z|)}{\varphi_{\nu,\alpha}(|z|)},$$
with equality when $z=|z|$. The previous inequality along with the minimum principle for harmonic functions and the fact that $\lim_{z\to0} \left.{z\varphi_{\nu,\alpha}'(z)}\right/{\varphi_{\nu,\alpha}(z)}=1$ complete the proof of the theorem.
\end{proof}
	
\begin{proof}[\bf Proof of Theorem \ref{Theorem3}]
Taking the logarithm of both sides of \eqref{Nor_CWF2} and differentiating with respect to $z$, we obtain
\begin{equation}\label{log_der_g}
\frac{g_{L,\eta}'(z)}{g_{L,\eta}(z)}=-\frac{L}{z}+\frac{F_{L,\eta}'(z)}{F_{L,\eta}(z)}.
\end{equation}
Multiplying both sides by $ze^{i\theta}$ and comparing the real parts we find that
\begin{equation}\label{spiral_g}
\real\left(e^{i\theta}\frac{zg_{L,\eta}'(z)}{g_{L,\eta}(z)}\right)=-\real\left(e^{i\theta}L\right)+\real\left(e^{i\theta}\frac{zF_{L,\eta}'(z)}{F_{L,\eta}(z)}\right).
\end{equation}
For the starlikeness property we substitute $\theta=0$ in \eqref{spiral_g} to obtain
\begin{equation}\label{real_log_der_g}
\real\left(\frac{zg_{L,\eta}'(z)}{g_{L,\eta}(z)}\right)=-\real L+\real\left(\frac{zF_{L,\eta}'(z)}{F_{L,\eta}(z)}\right).
\end{equation}
By using \eqref{final_reln} and \eqref{real_log_der_g} it turns out to be
\begin{equation}\label{rl_logder_g}
\real\left(\frac{zg_{L,\eta}'(z)}{g_{L,\eta}(z)}\right)\geq-\real L+\frac{|z|F_{l,\eta}'(|z|)}{F_{l,\eta}(|z|)}
\end{equation}
for all $|z|<r$. For $L=x+\mathrm{i}y$, the condition $l(l+1)=\real\left[L(L+1)\right]$ where $l>-1$ implies $\real\left[L(L+1)\right]>-\frac{1}{4}$ which is restated as $y^2<x(x+1)+\frac{1}{4}$. Moreover, \eqref{condn} implies that the right-hand side of \eqref{rl_logder_g} is positive for sufficiently small values of $|z|$. From this we conclude that
\begin{align*}
\real\left(\frac{zg_{L,\eta}'(z)}{g_{L,\eta}(z)}\right)\geq 0,\qquad |z|\leq\rho^*_{l,\eta,1},
\end{align*}
where $\rho^*_{l,\eta,1}$ is the smallest positive zero of the function
$$r\mapsto rF_{l,\eta}'(r)-\left(\real L\right)F_{l,\eta}(r).$$
For non-negative real values of $L$ the fact that $g_{L,\eta}'(z)$ vanishes for $\rho^*_{l,\eta,1}$ implies that $z\mapsto g_{L,\eta}(z)$ is not univalent in any domain like $|z|=\rho>\rho^*_{l,\eta,1}$. Hence, the normalized Coulomb wave function $z\mapsto g_{L,\eta}(z)$ is regular, univalent and spirallike in disk $|z|<\rho^*_{l,\eta,1}$.

If $L>-1,$ then \eqref{log_der_g} and \eqref{rl_logder_g} imply that
$$\real\left(\frac{zg_{L,\eta}'(z)}{g_{L,\eta}(z)}\right)\geq-L+\real\left(\frac{|z|F_{L,\eta}'(|z|)}{F_{L,\eta}(|z|)}\right)=\real\left(\frac{|z|g_{L,\eta}'(|z|)}{g_{L,\eta}(|z|)}\right).$$
The above inequality together with the fact that $\lim_{z\to0} \left.{zg_{L,\eta}'(z)}\right/{g_{L,\eta}(z)}=1$ and the minimum principle for harmonic functions imply that
$$\real\left(\frac{zg_{L,\eta}'(z)}{g_{L,\eta}(z)}\right)>\beta$$
holds if and only if $|z|<\tau^*_{l,\eta,1}$ where $\tau^*_{l,\eta,1}$ is the smallest positive root of the equation
$$-\real L+\frac{rF_{L,\eta}'(r)}{F_{L,\eta}(r)}=\beta,$$
which is equivalent to $$rF_{L,\eta}'(r)-\left(\beta+\real L\right)F_{L,\eta}(r)=0.$$
\end{proof}

\begin{proof}[\bf Proof of Lemma \ref{lemma_asy1}]
By using mathematical induction on $k$, we establish that for positive integer $k$ and real $L>k+1$, the Rayleigh sum $Z_\eta^{\left(k\right)}\left(L\right)$, defined in equation \eqref{Z_Coul}, can be expressed in the form \eqref{z_even_series} and \eqref{z_odd_series}. To do this we use the following recurrence relation presented in the work of \v{S}tampach and \v{S}\'{t}ov\'{\i}\v{c}ek \cite[eq. 79]{Flen_OrthPol_CWF}
$$Z_\eta^{\left(k+1\right)}\left(L\right)=\frac{1}{2L+k+2}\left(\frac{2\eta}{L+1}Z_\eta^{\left(k\right)}\left(L\right)+\sum_{l=1}^{k-2}Z_\eta^{\left(l+1\right)}\left(L\right)Z_\eta^{\left(k-l\right)}\left(L\right)\right),$$
where $k\in\{2,3,\ldots\}$ and the particular case \cite[eq. 78]{Flen_OrthPol_CWF}
\begin{equation}\label{ev_zta}
		Z_\eta^{\left(2\right)}\left(L\right)=\frac{1}{2L+3}\left(1+\frac{\eta^2}{\left(L+1\right)^2}\right).
\end{equation}
To deal with even and odd values of $k$ we write the above recurrence relation as
\begin{equation}\label{evn_zta}
Z_\eta^{\left(2k\right)}\left(L\right)=\frac{1}{2L+2k+1}\left(\frac{2\eta}{L+1}Z_\eta^{\left(2k-1\right)}\left(L\right)+\sum_{l=1}^{2k-3}Z_\eta^{\left(l+1\right)}\left(L\right)Z_\eta^{\left(2k-l-1\right)}
\left(L\right)\right), \quad k\in\{2,3,\ldots\},
\end{equation}
\begin{equation}\label{odd_zta}
Z_\eta^{\left(2k+1\right)}\left(L\right)=\frac{1}{2L+2k+2}\left(\frac{2\eta}{L+1}Z_\eta^{\left(2k\right)}\left(L\right)+\sum_{l=1}^{2k-2}Z_\eta^{\left(l+1\right)}\left(L\right)Z_\eta^{\left(2k-l\right)}\left(L\right)\right), \quad k\in \mathbb{N}.
\end{equation}
Now, for $L>\frac{\alpha+1}{2}$ we obtain that
\begin{equation}\label{p_n_series}
\frac{1}{2L+\alpha+1}=\frac{1}{L}\sum_{n\geq0}\frac{p_n^{\left(\alpha\right)}}{L^n},
\end{equation}
 where
$p_n^{\left(\alpha\right)}$ is given by \eqref{p_n}. Expanding the right-hand side of \eqref{ev_zta} and using \eqref{p_n_series} we arrive at
\begin{align*}
Z_\eta^{\left(2\right)}\left(L\right)&=\frac{1}{L}\sum_{n\geq0}\frac{p_n^{\left(2\right)}}{L^n}\left(1+\frac{\eta^2}{L^2\left(1+1/L\right)^2}\right)\\
&=\frac{1}{L}\sum_{n\geq0}\frac{p_n^{\left(2\right)}}{L^n}\left(1+\sum_{n\geq0}\left(-1\right)^n\frac{(n+1)\eta^2}{L^{n+2}}\right)\\
&=\frac{1}{L}\left(\sum_{n\geq0}\frac{p_n^{\left(2\right)}}{L^n}+\sum_{n\geq0}\sum_{m=0}^{n}{\left(-1\right)}^m\frac{\eta^2(m+1)p_{n-m}^{\left(2\right)}}{L^{n+2}}\right)\\
&=\frac{1}{L}\left(p_0^{\left(2\right)}+\frac{p_1^{\left(2\right)}}{L}+\sum_{n\geq0}\frac{p_{n+2}^{\left(2\right)}}{L^{n+2}}+\sum_{n\geq0}\sum_{m=0}^{n}{\left(-1\right)}^m\frac{\eta^2(m+1)p_{n-m}^{\left(2\right)}}{L^{n+2}}\right)\\
&=\frac{1}{L}\left(p_0^{\left(2\right)}+\frac{p_1^{\left(2\right)}}{L}+\sum_{n\geq0}\frac{1}{L^{n+2}}\left(p_{n+2}^{\left(2\right)}+\sum_{m=0}^{n}{\left(-1\right)}^m\eta^2(m+1)p_{n-m}^{\left(2\right)}\right)\right)\\
&=\frac{1}{L}\sum_{n\geq0}\frac{\zeta_{n,\eta}^{\left(2\right)}}{L^n},
\end{align*}
where $\zeta_{0,\eta}^{\left(2\right)}=p_0^{\left(2\right)}, \zeta_{1,\eta}^{\left(2\right)}=p_1^{\left(2\right)}$ and $$\zeta_{n+2,\eta}^{\left(2\right)}=p_{n+2}^{\left(2\right)}+\sum_{m=0}^{n}{\left(-1\right)}^m\eta^2(m+1)p_{n-m}^{\left(2\right)}$$ for $n\in \mathbb{N}_0$. By using \eqref{odd_zta} for $k=1$ and \eqref{p_n_series} we obtain that
\begin{align*}
Z_\eta^{\left(3\right)}\left(L\right)&=\frac{2\eta}{\left(2L+4\right)\left(L+1\right)}Z_\eta^{\left(2\right)}\left(L\right)=\frac{2\eta}{L^3}\left(\sum_{n\geq0}\frac{p_n^{\left(3\right)}}{L^n}\right)
\left(\sum_{n\geq0}\frac{\left(-1\right)^n}{L^n}\right)\left(\sum_{n\geq0}\frac{\zeta_{n,\eta}^{\left(2\right)}}{L^n}\right)\\
&=\frac{2\eta}{L^3}\left(\sum_{n\geq0}\sum_{l=0}^{n}\sum_{m=0}^{l}\frac{\left(-1\right)^m p_{l-m}^{\left(3\right)}\zeta_{n-l,\eta}^{\left(2\right)}}{L^n}\right)=\frac{1}{L^3}\sum_{n\geq0}\frac{\zeta_{n,\eta}^{\left(3\right)}}{L^n},
\end{align*}
where $$\zeta_{n,\eta}^{\left(3\right)}=2\eta \displaystyle{\sum_{l=0}^{n}\sum_{m=0}^{l}} \left(-1\right)^m p_{l-m}^{\left(3\right)}\zeta_{n-l,\eta}^{\left(2\right)}$$ for $n\in \mathbb{N}_0$.
Let $N\geq2$ and suppose that $Z_\eta^{\left(2k\right)}\left(L\right)$ and $Z_\eta^{\left(2k+1\right)}\left(L\right)$ can be expressed in the form \eqref{z_even_series} and \eqref{z_odd_series} for $1<k\leq N-1$, respectively. By using the expansions of $Z_\eta^{\left(2\right)}\left(L\right)$, $Z_\eta^{\left(3\right)}\left(L\right)$ and the induction hypothesis, we write \eqref{evn_zta} for $k=N$ as
\begin{align*}
Z_\eta^{\left(2N\right)}\left(L\right)&=\frac{2\eta}{\left(2L+2N+1\right)\left(L+1\right)}Z_\eta^{\left(2N-1\right)}\left(L\right)+\frac{1}{\left(2L+2N+1\right)}\sum_{l=1}^{2N-3}Z_\eta^{\left(l+1\right)}Z_\eta^{\left(2N-l-1\right)}\left(L\right).
\end{align*}
In view of \eqref{p_n_series} and the induction hypothesis, the first term in the right-hand side of the above equation can be expressed as
\begin{align*}
\frac{2\eta}{\left(2L+2N+1\right)\left(L+1\right)}Z_\eta^{\left(2N-1\right)}\left(L\right)&=\frac{2\eta}{L^{2N+1}}\left(\sum_{n\geq0}\frac{p_n^{\left(2N\right)}}{L^n}\right)
\left(\sum_{n\geq0}\frac{\left(-1\right)^n}{L^n}\right)\left(\sum_{n\geq0}\frac{\zeta_{n,\eta}^{\left(2N-1\right)}}{L^n}\right)\\
&=\frac{2\eta}{L^{2N+1}}\sum_{n\geq0}\sum_{l=0}^{n}\sum_{m=0}^{l}\frac{\left(-1\right)^mp_{l-m}^{\left(2N\right)}\zeta_{n-l,\eta}^{\left(2N-1\right)}}{L^n}\\
&=\frac{1}{L^{2N-1}}\sum_{n\geq0}\left(\frac{2\eta}{L^{n+2}}\sum_{l=0}^{n}\sum_{m=0}^{l}\left(-1\right)^mp_{l-m}^{\left(2N\right)}\zeta_{n-l,\eta}^{\left(2N-1\right)}\right),
\end{align*}
while the second term can be written as
\begin{align*}
&\frac{1}{\left(2L+2N+1\right)}\sum_{l=1}^{2N-3}Z_\eta^{\left(l+1\right)}Z_\eta^{\left(2N-l-1\right)}\left(L\right)\\
&=\frac{1}{2L+2N+1}\left(\sum_{l=1}^{N-2}Z_\eta^{\left(2l+1\right)}\left(L\right)Z_\eta^{\left(2N-2l-1\right)}\left(L\right)+\sum_{l=0}^{N-2}Z_\eta^{\left(2l+2\right)}\left(L\right)Z_\eta^{\left(2N-2l-2\right)}\left(L\right)\right)\\
&=\frac{1}{L}\sum_{n\geq0}\frac{p_n^{\left(2N\right)}}{L^n}\left[\sum_{l=1}^{N-2}\frac{1}{L^{2N}}\left(\sum_{n\geq0}\frac{\zeta_{n,\eta}^{\left(2l+1\right)}}{L^n}\right)\left(\sum_{m\geq0}\frac{\zeta_{m,\eta}^{\left(2N-2l-1\right)}}{L^m}\right)\right.\\
&\left.\quad\quad\quad+\sum_{l=0}^{N-2}\frac{1}{L^{2N-2}}\left(\sum_{n\geq0}\frac{\zeta_{n,\eta}^{\left(2l+2\right)}}{L^n}\right)\left(\sum_{n\geq0}\frac{\zeta_{n,\eta}^{\left(2N-2l-2\right)}}{L^n}\right)\right]\\
&=\frac{1}{L^{2N-1}}\left[\frac{1}{L^2}\left(\sum_{n\geq0}\frac{p_n^{(2N)}}{L^n}\right)\sum_{l=1}^{N-2}\left(\sum_{n\geq0}\sum_{m=0}^{n}\frac{\zeta_{m,\eta}^{\left(2l+1\right)}\zeta_{n-m,\eta}^{\left(2N-2l-1\right)}}{L^n}\right)\right.\\     &\left.\quad\quad\quad+\left(\sum_{n\geq0}\frac{p_n^{(2N)}}{L^n}\right)\sum_{l=0}^{N-2}\left(\sum_{n\geq0}\sum_{m=0}^{n}\frac{\zeta_{m,\eta}^{\left(2l+2\right)}\zeta_{n-m,\eta}^{\left(2N-2l-2\right)}}{L^n}\right)\right]
\end{align*}
or equivalently
\begin{align*}
&\frac{1}{L^{2N-1}}\left[\frac{1}{L^2}\sum_{l=1}^{N-2}\left(\sum_{n\geq0}\sum_{q=0}^{n}\sum_{m=0}^{q}\frac{\zeta_{m,\eta}^{\left(2l+1\right)}\zeta_{q-m,\eta}^{\left(2N-2l-1\right)}p_{n-q}^{(2N)}}{L^n}\right)\right.\\     &\left.\quad\quad\quad+\sum_{l=0}^{N-2}\left(\sum_{n\geq0}\sum_{q=0}^{n}\sum_{m=0}^{q}\frac{\zeta_{m,\eta}^{\left(2l+2\right)}\zeta_{q-m,\eta}^{\left(2N-2l-2\right)}p_{n-q}^{(2N)}}{L^n}\right) \right]\\
&=\frac{1}{L^{2N-1}}\left[\sum_{l=1}^{N-2}\left(\sum_{n=0}^{\infty}\sum_{q=0}^{n}\sum_{m=0}^{q}\frac{\zeta_{m,\eta}^{\left(2l+1\right)}\zeta_{q-m,\eta}^{\left(2N-2l-1\right)}p_{n-q}^{(2N)}}{L^{n+2}}\right)  +\sum_{l=0}^{N-2}\zeta_{0,\eta}^{\left(2l+2\right)}\zeta_{0,\eta}^{\left(2N-2l-2\right)}p_0^{(2N)}\right.\\ &\left.\quad+\sum_{l=0}^{N-2}\sum_{q=0}^{1}\sum_{m=0}^{q}\frac{\zeta_{m,\eta}^{\left(2l+2\right)}\zeta_{q-m,\eta}^{\left(2N-2l-2\right)}p_{1-q}^{(2N)}}{L} +\sum_{l=0}^{N-2}\left(\sum_{n\geq0}\sum_{q=0}^{n+2}\sum_{m=0}^{q}\frac{\zeta_{m,\eta}^{\left(2l+2\right)}\zeta_{q-m,\eta}^{\left(2N-2l-2\right)}p_{n+2-q}^{(2N)}}{L^{n+2}}\right) \right]\\
&=\frac{1}{L^{2N-1}}\left[\sum_{l=0}^{N-2}\zeta_{0,\eta}^{\left(2l+2\right)}\zeta_{0,\eta}^{\left(2N-2l-2\right)}p_0^{(2N)}+\sum_{l=0}^{N-2}\sum_{q=0}^{1}\sum_{m=0}^{q}\frac{\zeta_{m,\eta}^{\left(2l+2\right)}
\zeta_{q-m,\eta}^{\left(2N-2l-2\right)}p_{1-q}^{(2N)}}{L}\right.\\
&\left.\quad+\sum_{n\geq0}\frac{1}{L^{n+2}}\left(\sum_{l=1}^{N-2}\sum_{q=0}^{n}\sum_{m=0}^{q}\zeta_{m,\eta}^{\left(2l+1\right)}\zeta_{q-m,\eta}^{\left(2N-2l-1\right)}p_{n-q}^{(2N)} +\sum_{l=0}^{N-2}\sum_{q=0}^{n+2}\sum_{m=0}^{q}\zeta_{m,\eta}^{\left(2l+2\right)}\zeta_{q-m,\eta}^{\left(2N-2l-2\right)}p_{n+2-q}^{(2N)}\right)\right].
\end{align*}
Combining the above expressions we obtain that
\begin{align*}
Z_\eta^{\left(2N\right)}&\left(L\right)=\frac{1}{L^{2N-1}}\left[\sum_{l=0}^{N-2}\zeta_{0,\eta}^{\left(2l+2\right)}\zeta_{0,\eta}^{\left(2N-2l-2\right)}p_0^{(2N)}+
\sum_{l=0}^{N-2}\sum_{q=0}^{1}\sum_{m=0}^{q}\frac{\zeta_{m,\eta}^{\left(2l+2\right)}\zeta_{q-m,\eta}^{\left(2N-2l-2\right)}p_{1-q}^{(2N)}}{L}\right.\\
&+\sum_{n\geq0}\frac{1}{L^{n+2}}\left(\sum_{l=1}^{N-2}\sum_{q=0}^{n}\sum_{m=0}^{q}\zeta_{m,\eta}^{\left(2l+1\right)}\zeta_{q-m,\eta}^{\left(2N-2l-1\right)}p_{n-q}^{(2N)}
+\sum_{l=0}^{N-2}\sum_{q=0}^{n+2}\sum_{m=0}^{q}\zeta_{m,\eta}^{\left(2l+2\right)}\zeta_{q-m,\eta}^{\left(2N-2l-2\right)}p_{n+2-q}^{(2N)}\right.\\
&\left.\left.+2\eta\sum_{l=0}^{n}\sum_{m=0}^{l}\left(-1\right)^mp_{l-m}^{\left(2N\right)}\zeta_{n-l,\eta}^{\left(2N-1\right)}\right)\right]
\end{align*}
or equivalently we arrive at
\begin{equation}\label{Z_2N_series}
Z_\eta^{\left(2N\right)}\left(L\right)=\frac{1}{L^{2N-1}}\sum_{n\geq0}\frac{\zeta_{n,\eta}^{\left(2N\right)}}{L^n},
\end{equation}
where
$$\zeta_{0,\eta}^{\left(2N\right)}=\sum_{l=0}^{N-2}\zeta_{0,\eta}^{\left(2l+2\right)}\zeta_{0,\eta}^{\left(2N-2l-2\right)}p_0^{(2N)},$$
$$\zeta_{1,\eta}^{\left(2N\right)}=\sum_{l=0}^{N-2}\sum_{q=0}^{1}\sum_{m=0}^{q}\zeta_{m,\eta}^{\left(2l+2\right)}\zeta_{q-m,\eta}^{\left(2N-2l-2\right)}p_{1-q}^{(2N)}$$
and
\begin{align*}
\zeta_{n+2,\eta}^{\left(2N\right)}&=\sum_{l=1}^{N-2}\sum_{q=0}^{n}\sum_{m=0}^{q}\zeta_{m,\eta}^{\left(2l+1\right)}\zeta_{q-m,\eta}^{\left(2N-2l-1\right)}p_{n-q}^{(2N)} +\sum_{l=0}^{N-2}\sum_{q=0}^{n+2}\sum_{m=0}^{q}\zeta_{m,\eta}^{\left(2l+2\right)}\zeta_{q-m,\eta}^{\left(2N-2l-2\right)}p_{n+2-q}^{(2N)}\\
&\quad\quad+2\eta\sum_{l=0}^{n}\sum_{m=0}^{l}\left(-1\right)^mp_{l-m}^{\left(2N\right)}\zeta_{n-l,\eta}^{\left(2N-1\right)}
\end{align*}
for $n\in \mathbb{N}_0.$ By using $Z_\eta^{\left(2\right)}\left(L\right)$, $Z_\eta^{\left(3\right)}\left(L\right)$ and the induction hypothesis, we write \eqref{odd_zta} for $k=N$ as
$$Z_\eta^{\left(2N+1\right)}\left(L\right)=\frac{2\eta}{\left(2L+2N+2\right)\left(L+1\right)}Z_\eta^{\left(2N\right)}\left(L\right)+\frac{1}{\left(2L+2N+2\right)}
\sum_{l=1}^{2N-2}Z_\eta^{\left(l+1\right)}\left(L\right)Z_\eta^{\left(2N-l\right)}\left(L\right).$$
In view of \eqref{p_n_series} and the induction hypothesis again, the first term in the right-hand side of above equation can be expressed as
\begin{align*}
\frac{2\eta}{\left(2L+2N+2\right)\left(L+1\right)}Z_\eta^{\left(2N\right)}\left(L\right)&=\frac{2\eta}{L^{2N+1}}\left(\sum_{n\geq0}\frac{p_n^{\left(2N+1\right)}}{L^n}\right)
\left(\sum_{n\geq0}\frac{\left(-1\right)^n}{L^n}\right)\left(\sum_{n\geq0}\frac{\zeta_{n,\eta}^{\left(2N\right)}}{L^n}\right)\\
&=\frac{1}{L^{2N+1}}\sum_{n\geq0}\left(\frac{2\eta}{L^{n}}\sum_{l=0}^{n}\sum_{m=0}^{l}\left(-1\right)^mp_{l-m}^{\left(2N+1\right)}\zeta_{n-l,\eta}^{\left(2N\right)}\right),
\end{align*}
while the second term can be written as
\begin{align*}
&\frac{1}{\left(2L+2N+2\right)}\sum_{l=1}^{2N-2}Z_\eta^{\left(l+1\right)}\left(L\right)Z_\eta^{\left(2N-l\right)}\left(L\right)\\
&=\frac{1}{2L+2N+2}\left(\sum_{l=1}^{N-1}Z_\eta^{\left(2l+1\right)}\left(L\right)Z_\eta^{\left(2N-2l\right)}\left(L\right)
+\sum_{l=0}^{N-2}Z_\eta^{\left(2l+2\right)}\left(L\right)Z_\eta^{\left(2N-2l-1\right)}\left(L\right)\right)\\
&=\left(\frac{1}{L}\sum_{n\geq0}\frac{p_n^{\left(2N+1\right)}}{L^n}\right)\left[\sum_{l=1}^{N-1}\frac{1}{L^{2N}}\left(\sum_{n\geq0}\frac{\zeta_{n,\eta}^{\left(2l+1\right)}}{L^n}\right)
\left(\sum_{m\geq0}\frac{\zeta_{m,\eta}^{\left(2N-2l\right)}}{L^m}\right)\right.\\
&\left.\quad\quad+\sum_{l=0}^{N-2}\frac{1}{L^{2N}}\left(\sum_{n\geq0}\frac{\zeta_{n,\eta}^{\left(2l+2\right)}}{L^n}\right)
\left(\sum_{n\geq0}\frac{\zeta_{n,\eta}^{\left(2N-2l-1\right)}}{L^n}\right)\right]\\
&=\frac{1}{L^{2N+1}}\left[\left(\sum_{n\geq0}\frac{p_n^{\left(2N+1\right)}}{L^n}\right)\sum_{l=1}^{N-1}
\left(\sum_{n\geq0}\sum_{m=0}^{n}\frac{\zeta_{m,\eta}^{\left(2l+1\right)}\zeta_{n-m,\eta}^{\left(2N-2l\right)}}{L^n}\right)\right.\\
&\left.\quad\quad+\left(\sum_{n\geq0}\frac{p_n^{\left(2N+1\right)}}{L^n}\right)\sum_{l=0}^{N-2}\left(\sum_{n\geq0}\sum_{m=0}^{n}
\frac{\zeta_{m,\eta}^{\left(2l+2\right)}\zeta_{n-m,\eta}^{\left(2N-2l-1\right)}}{L^n}\right)\right]\\
&=\frac{1}{L^{2N+1}}\left[\sum_{n\geq0}\frac{1}{L^n}\left(\sum_{l=1}^{N-1}\sum_{q=0}^{n}\sum_{m=0}^{q}\zeta_{m,\eta}^{\left(2l+1\right)}\zeta_{q-m,\eta}^{\left(2N-2l\right)}p_{n-q}^{\left(2N+1\right)} +\sum_{l=0}^{N-2}\sum_{q=0}^{n}\sum_{m=0}^{q}\zeta_{m,\eta}^{\left(2l+2\right)}\zeta_{q-m,\eta}^{\left(2N-2l-1\right)}p_{n-q}^{2N+1}\right)\right].
\end{align*}
Combining the above expressions we write
\begin{align*}
Z_\eta^{\left(2N+1\right)}\left(L\right)=\frac{1}{L^{2N+1}}&\left[\sum_{n\geq0}\frac{1}{L^n}\left(2\eta\sum_{l=0}^{n}\sum_{m=0}^{l}\left(-1\right)^mp_{l-m}^{\left(2N+1\right)}\zeta_{n-l,\eta}^{\left(2N\right)}	
	+\sum_{l=1}^{N-1}\sum_{q=0}^{n}\sum_{m=0}^{q}\zeta_{m,\eta}^{\left(2l+1\right)}\zeta_{q-m,\eta}^{\left(2N-2l\right)}p_{n-q}^{\left(2N+1\right)}\right.\right.\\ &\left.\left.\quad\quad+\sum_{l=0}^{N-2}\sum_{q=0}^{n}\sum_{m=0}^{q}\zeta_{m,\eta}^{\left(2l+2\right)}\zeta_{q-m,\eta}^{\left(2N-2l-1\right)}p_{n-q}^{2N+1}\right)\right],
\end{align*}
which can be written as
\begin{equation}\label{Z_2N+1_series}
Z_\eta^{\left(2N+1\right)}\left(L\right)=\frac{1}{L^{2N+1}}\sum_{n\geq0}\frac{\zeta_{n,\eta}^{\left(2N+1\right)}}{L^n},
\end{equation}
where
\begin{align*}
\zeta_{n,\eta}^{\left(2N+1\right)}=2\eta&\sum_{l=0}^{n}\sum_{m=0}^{l}\left(-1\right)^mp_{l-m}^{\left(2N+1\right)}\zeta_{n-l,\eta}^{\left(2N\right)}+\sum_{l=1}^{N-1}\sum_{q=0}^{n}
\sum_{m=0}^{q}\zeta_{m,\eta}^{\left(2l+1\right)}\zeta_{q-m,\eta}^{\left(2N-2l\right)}p_{n-q}^{\left(2N+1\right)}\\ &\quad\quad+\sum_{l=0}^{N-2}\sum_{q=0}^{n}\sum_{m=0}^{q}\zeta_{m,\eta}^{\left(2l+2\right)}\zeta_{q-m,\eta}^{\left(2N-2l-1\right)}p_{n-q}^{2N+1}
\end{align*}
for $n\in \mathbb{N}_0.$ Thus in view of equations \eqref{Z_2N_series} and \eqref{Z_2N+1_series} the proof is completed.
\end{proof}

\begin{proof}[\bf Proof of Lemma \ref{lemma_asy2}]
We shall use the next relation deduced by \v{S}tampach and \v{S}\'{t}ov\'{\i}\v{c}ek \cite[p. 253]{Flen_OrthPol_CWF}
\begin{align}\label{rec_der}
\frac{2\left(\rho-\eta\right)}{\rho^2-2\eta\rho-L\left(L+1\right)}&\left(\sum_{n\geq2}\tilde{Z}_\eta^{\left(n\right)}\left(L\right)\rho^{n+1}-L\rho-\frac{\left(L+2\right)\eta\rho^2}{\left(L+1\right)^2}\right)
+\left(\sum_{n\geq2}\tilde{Z}_\eta^{\left(n\right)}\left(L\right)\rho^n-L-\frac{\left(L+2\right)\eta\rho}{\left(L+1\right)^2}\right)^2\\
&=L^2+\frac{2\left(L^2+L-1\right)\eta\rho}{\left(L+1\right)^2}-\rho^2+\sum_{n\geq2}\left(n+1\right)\tilde{Z}_\eta^{\left(n\right)}\left(L\right)\rho^n.\nonumber
\end{align}
and for the sake of convenience we write
$$\frac{2\left(\rho-\eta\right)}{\rho^2-2\eta\rho-L\left(L+1\right)}=\sum_{n\geq0}a_n\rho^n.$$
Observe that the left-hand side of \eqref{rec_der} can be written as
\begin{align*}
-L\sum_{n\geq0}a_n&\rho^{n+1}-\frac{\left(L+2\right)\eta}{\left(L+1\right)^2}\sum_{n\geq0}a_n\rho^{n+2}+\sum_{n\geq0}a_n\rho^n
\sum_{n\geq0}\tilde{Z}_\eta^{\left(n+2\right)}\left(L\right)\rho^{n+3}+\left(L+\frac{\left(L+2\right)}{\left(L+1\right)^2}\eta\rho\right)^2\\
&+\left(\sum_{n\geq0}\tilde{Z}_\eta^{\left(n+2\right)}\rho^{n+2}\right)^2-2\left(L+\frac{\left(L+2\right)\eta}{\left(L+1\right)^2}\rho\right)
\left(\sum_{n\geq0}\tilde{Z}_\eta^{\left(n+2\right)}\left(L\right)\rho^{n+2}\right).
\end{align*}
Moreover, we use the notation $\frac{\left(L+2\right)}{\left(L+1\right)^2}\eta=\tilde{p}$. With this notation, the above expression reads as follows
\begin{align*}
-L&\sum_{n\geq0}a_n\rho^{n+1}-\tilde{p}\sum_{n\geq0}a_n\rho^{n+2}+\sum_{n\geq0}\sum_{m=0}^{n}a_m\tilde{Z}_\eta^{\left(n-m+2\right)}\left(L\right)\rho^{n+3}+\left(L^2+\tilde{p}^2\rho^2+2L\tilde{p}\rho\right)\\
&+\sum_{n\geq0}\sum_{m=0}^{n}\tilde{Z}_\eta^{\left(m+2\right)}\left(L\right)\tilde{Z}_\eta^{\left(n-m+2\right)}\left(L\right)\rho^{n+4}-2L\sum_{n\geq0}\tilde{Z}_\eta^{\left(n+2\right)}\left(L\right)\rho^{n+2}
-2\tilde{p}\sum_{n\geq0}\tilde{Z}_\eta^{\left(n+2\right)}\left(L\right)\rho^{n+3}.
\end{align*}
Replacing the left-hand side of \eqref{rec_der} by the above expression and equating the coefficients of $\rho$ we get the required recurrence relation for $\tilde{Z}_\eta^{\left(k\right)}\left(L\right)$.
In the current context, we formulate the so-called Euler-Rayleigh inequalities \cite[p. 252]{Flen_OrthPol_CWF} to establish the bounds for the radius of starlikeness associated with the Coulomb wave function when $\eta< 0$ as
\begin{align*}
\left(\tilde{Z}_\eta^{\left(2s\right)}\left(L\right)\right)^{-\frac{1}{s}}<\tilde{\rho}_{L,1}^2<\frac{\tilde{Z}_\eta^{\left(2s\right)}\left(L\right)}{\tilde{Z}_\eta^ {\left(2s+2\right)}\left(L\right)}.
\end{align*}
The aforementioned inequalities together with the observation that the first zero in modulus of the derivative of the Coulomb wave function coincides with the radius of starlikeness of the Coulomb wave function, constitutes the final element in the proof of the Lemma.
\end{proof}

\begin{proof}[\bf Proof of Lemma \ref{Lemma3}]
By using the fact that $p_0^{\left(\alpha\right)}=\frac{1}{2}$ in Lemma \ref{lemma_asy1}, we have that
$$\zeta_{0,\eta}^{\left(2k\right)}=\frac{1}{2}\sum_{l=0}^{k-2}\zeta_{0,\eta}^{\left(2l+2\right)}\zeta_{0,\eta}^{\left(2k-2l-2\right)}.$$
We use mathematical induction on $k$ to prove
\begin{equation}\label{const_evn_zeta}
\zeta_{0,\eta}^{\left(2k\right)}\leq\binom{2k}{k}\frac{1}{2^{2k}\left(2k-1\right)}.
\end{equation}	
For $k=1$, $\zeta_{0,\eta}^{\left(2\right)}=p_0^{\left(2\right)}=\frac{1}{2}=\frac{1}{4}\binom{2}{1}.$ Suppose that the claim is true for $1\leq k\leq N-1$. By using the induction hypothesis, we obtain for $k=N$
\begin{align*}
\zeta_{0,\eta}^{\left(2N\right)}&=\frac{1}{2}\sum_{l=0}^{N-2}\zeta_{0,\eta}^{\left(2l+2\right)}\zeta_{0,\eta}^{\left(2N-2l-2\right)}\\
&\leq\frac{1}{2}\sum_{l=0}^{N-2}\binom{2l+2}{l+1}\frac{1}{2^{2l+2}\left(2l+1\right)}\binom{2N-2l-2}{N-l-1}\frac{1}{2^{2N-2l-2}\left(2N-2l-3\right)}\\
&=\frac{1}{2^{2N+1}}\sum_{l=0}^{N-2}\binom{2l+2}{l+1}\frac{1}{\left(2l+1\right)}\binom{2N-2l-2}{N-l-1}\frac{1}{\left(2N-2l-3\right)}\\
&\leq\frac{1}{2^{2N+1}}\frac{1}{2N-1}\binom{2N}{N}<\frac{1}{2^{2N}}\frac{1}{2N-1}\binom{2N}{N}\leq\frac{1}{2N-1}<1,
\end{align*}
where we used the Chu-Vandermonde identity and for the last line inequality we refer to the proof of \cite[Lemma 2]{BN21}. Therefore by mathematical induction, we write $\zeta_{0,\eta}^{\left(2k\right)}=\mathcal{O}\left(1\right)$ for all $k\in\mathbb{N}$. By using the equation \eqref{z_even_series} for large $L$ we obtain \eqref{even_zeta}. Now, consider the zeta function associated with odd power of zeros of Coulomb wave functions, namely $Z_\eta^{\left(2k+1\right)}\left(L\right)$.
From Lemma \ref{lemma_asy1} we obtain that
$$\zeta_{0,\eta}^{\left(2k+1\right)}=2\eta p_{0}^{\left(2k+1\right)}\zeta_{0,\eta}^{\left(2k\right)}
+\sum_{l=1}^{k-1}\zeta_{0,\eta}^{\left(2l+1\right)}\zeta_{0,\eta}^{\left(2k-2l\right)}p_{0}^{\left(2k+1\right)}
+\sum_{l=0}^{k-2}\zeta_{0,\eta}^{\left(2l+2\right)}\zeta_{0,\eta}^{\left(2k-2l-1\right)}p_{0}^{\left(2k+1\right)}.$$
Since $p_0^{\left(\alpha\right)}=\frac{1}{2},$ we can write the above equation as
\begin{equation}\label{expression_1}
\zeta_{0,\eta}^{\left(2k+1\right)}-\eta \zeta_{0,\eta}^{\left(2k\right)}=\frac{1}{2}\left(\sum_{l=1}^{k-1}\zeta_{0,\eta}^{\left(2l+1\right)}\zeta_{0,\eta}^{\left(2k-2l\right)}
+\sum_{l=0}^{k-2}\zeta_{0,\eta}^{\left(2l+2\right)}\zeta_{0,\eta}^{\left(2k-2l-1\right)}\right).
\end{equation}
Now, for $\eta \leq 0$ by using mathematical induction on $k$, we prove the inequalities
\begin{equation}\label{claim_2}
\zeta_{0,\eta}^{\left(2k+1\right)}\leq\zeta_{0,\eta}^{\left(2k+1\right)}-\eta \zeta_{0,\eta}^{\left(2k\right)}\leq\frac{1}{2^{2k+1}\left(2k\right)}|\eta|\binom{2k+1}{k+1}.
\end{equation}
Note that for $\eta \leq 0$, the first inequality in \eqref{claim_2} is trivially true. For $k=1$ the right-hand side of \eqref{expression_1} becomes zero as being an empty sum and therefore
$$\zeta_{0,\eta}^{\left(3\right)}-\eta\zeta_{0,\eta}^{\left(2\right)}=0.$$
Therefore the second inequality in \eqref{claim_2} is satisfied for $k=1$. Now, suppose that the second inequality in \eqref{claim_2} is true for $1\leq k\leq N-1$.
From equations \eqref{const_evn_zeta}, \eqref{expression_1} and by using the induction hypothesis, we obtain for $k=N$
\begin{align*}
\zeta_{0,\eta}^{\left(2N+1\right)}-\eta \zeta_{0,\eta}^{\left(2N\right)}&=\frac{1}{2}\left(\sum_{l=1}^{N-1}\zeta_{0,\eta}^{\left(2l+1\right)}\zeta_{0,\eta}^{\left(2N-2l\right)}+\sum_{l=0}^{N-2}\zeta_{0,\eta}^{\left(2l+2\right)}\zeta_{0,\eta}^{\left(2N-2l-1\right)}\right)\\
&\leq\frac{|\eta|}{2}\left(\sum_{l=1}^{N-1}\binom{2l+1}{l+1}\frac{1}{2^{2l+1}\left(2l\right)}\binom{2N-2l}{N-l}\frac{1}{2^{2N-2l}\left(2N-2l-1\right)}\right.\\
&\left.\qquad+\sum_{l=0}^{N-2}\binom{2l+2}{l+1}\frac{1}{2^{2l+2}\left(2l+1\right)}\binom{2N-2l-1}{N-l}\frac{1}{2^{2N-2l-1}\left(2N-2l-2\right)}\right)\\
&\leq\frac{|\eta|}{2}\left(\binom{2N+1}{N+1}\frac{1}{2^{2N+1}\left(2N\right)}+\binom{2N+1}{N+1}\frac{1}{2^{2N+1}\left(2N\right)}\right)\\
&=|\eta|\binom{2N+1}{N+1}\frac{1}{2^{2N+1}\left(2N\right)}.
\end{align*}
This completes the proof of the claim in \eqref{claim_2}. Now it is easy to verify that $\zeta_{0,\eta}^{2k+1}<\frac{9}{8}|\eta|$. We can write $\zeta_{0,\eta}^{\left(2k+1\right)}=\mathcal{O}\left(1\right)$ for all $k\in\mathbb{N}$. By using the equation \eqref{z_odd_series} for large $L$ we obtain \eqref{odd_zeta}.
\end{proof}

\begin{proof}[\bf Proof of Lemma \ref{Lemma4}]
Note that the coefficient $a_i$ in equation \eqref{coeff_a_n} can be expressed in the series form as given below
\begin{align*}
a_{2n}=\frac{1}{L^{2n+2}}\sum_{m\geq0}\frac{a_m^{(2n)}}{L^m}\quad \mbox{and} \quad a_{2n+1}=\frac{1}{L^{2n+2}}\sum_{m\geq0}\frac{a_m^{(2n+1)}}{L^m}\quad \mbox{for all} \quad n\in \mathbb{N}_0.
\end{align*}
Moreover, we can write $\tilde{p}=\frac{(L+2)}{(L+1)^2}\eta=\frac{1}{L}\sum_{n\geq0}\frac{p_n}{L^n}$ for some coefficient $p_n$. By using these expressions, analogous to Lemma \ref{lemma_asy1}, we can express
$$\tilde{Z}^{\left(2k\right)}\left(L\right)=\frac{1}{L^{2k-1}}\sum_{n\geq0}\frac{\tilde{\zeta}_{n,\eta}^{\left(2k\right)}}{L^n}$$
and
$$\tilde{Z}^{\left(2k+1\right)}\left(L\right)=\frac{1}{L^{2k+1}}\sum_{m\geq0}\frac{\tilde{\zeta}_{m,\eta}^{\left(2k+1\right)}}{L^m}.$$
For arbitrary $s\geq1$ we have
$$\frac{\tilde{Z}_\eta^{\left(2s\right)}\left(L\right)}{\tilde{Z}_\eta^{\left(2s+2\right)}\left(L\right)}
=\frac{1}{L^{2s-1}}\sum_{n\geq0}\frac{\tilde{\zeta}_{n,\eta}^{\left(2s\right)}}{L^n}\left/\frac{1}{L^{2s+1}}
\sum_{m\geq0}\frac{\tilde{\zeta}_{m,\eta}^{\left(2s+2\right)}}{L^m}=L^2\sum_{n\geq0}\frac{b_n}{L^n}\right.$$
and thus we can write $$\frac{\tilde{Z}_\eta^{\left(2s\right)}\left(L\right)}{\tilde{Z}_\eta^{\left(2s+2\right)}\left(L\right)}=L^2\left(c'+\mathcal{O}\left(\frac{1}{L}\right)\right),$$
where $c'$ is some positive constant. Since the inequality \eqref{ros_ineqn} becomes equality as $s\to \infty$, from the above equation and \eqref{ros_ineqn}, we obtain  $$\left(r^*\left(f_{L,\eta}\right)\right)^2=L^2\left(c'+\mathcal{O}\left(\frac{1}{L}\right)\right).$$
Consequently we can express the radius of starlikeness as follows $$r^*\left(f_{L,\eta}\right)=L\left(c+\mathcal{O}\left(\frac{1}{L}\right)\right),$$
where $c$ is some positive constant. It is worth to note that we will get similar asymptotic form for $r^*\left(f_{L,\eta}\right)$ when we consider the left-side of inequality \eqref{ros_ineqn} for $s\to \infty$ and large $L$.
\end{proof}

\begin{proof}[\bf Proof of Theorem \ref{asym_the}]
The Weierstrass canonical product expansion of the Coulomb wave function reads as (see \cite[eqn. 76]{Flen_OrthPol_CWF})
$$F_{L,\eta}\left(z\right)=C_L\left(\eta\right)z^{\left(L+1\right)}e^{\frac{\eta z}{L+1}}{\displaystyle \prod_{n\geq1}^{}\left(1-\frac{z}{\rho_{L,n}}\right)e^{\frac{z}{\rho_{L,n}}}},$$
where $\rho_{L,n}$ is the $n$th zero of the Coulomb wave function $F_{L,\eta}(z)$. 	
Taking the logarithm of both sides of the above equation and then differentiating with respect to $z$, we conclude that
\begin{equation}\label{log_der3}
\frac{f_{L,\eta}^\prime\left(z\right)}{f_{L,\eta}\left(z\right)}=\frac{1}{L+1}\frac{F_{L,\eta}^\prime\left(z\right)}{F_{L,\eta}\left(z\right)}
=\frac{1}{z}+\frac{\eta}{\left(L+1\right)^2}-\frac{1}{\left(L+1\right)}\sum_{n\geq1}\frac{z}{\rho_{L,n}\left(\rho_{L,n}-z\right)}.
\end{equation}
Since the radius of starlikeness of the normalized Coulomb wave function $z\mapsto f_{L,\eta}(z)$ is the smallest zero in modulus of $z\mapsto F_{L,\eta}^\prime(z)$, the left-hand side of equation \eqref{log_der3} vanishes at $r^*\left(f_{L,\eta}\right)$ and thus we have
\begin{equation}\label{eqn23}
\frac{1}{r^*\left(f_{L,\eta}\right)}=\frac{1}{L+1}\sum_{n\geq1}\frac{r^*\left(f_{L,\eta}\right)}{\rho_{L,n}\left(\rho_{L,n}
-r^*\left(f_{L,\eta}\right)\right)}-\frac{\eta}{\left(L+1\right)^2}.
\end{equation}
Now with the help of Lemma \ref{Lemma4} we can write
$$r^*\left(f_{L,\eta}\right)=L\left(c+\mathcal{O}\left(\frac{1}{L}\right)\right) = L\left(c+\epsilon\left(L\right)\right)$$
for large $L$, where $c$ is some constant and $\epsilon\left(L\right)=\mathcal{O}\left(\frac{1}{L}\right)$.
Rearranging \eqref{eqn23} we obtain that
\begin{align*}
\frac{1}{L\left(c+\epsilon\left(L\right)\right)}&=\frac{1}{L+1}\sum_{n\geq1}\frac{L\left(c+\epsilon\left(L\right)\right)}{\rho_{L,n}\left(\rho_{L,n}-L\left(c+\epsilon\left(L\right)\right)\right)}
-\frac{\eta}{\left(L+1\right)^2}\\&=\frac{1}{L+1}\sum_{n\geq1}\frac{L\left(c+\epsilon\left(L\right)\right)}{\rho_{L,n}^2\left(1-\frac{L\left(c+\epsilon\left(L\right)\right)}{\rho_{L,n}}\right)}
-\frac{\eta}{\left(L+1\right)^2}\\
&=\frac{1}{L+1}\sum_{n\geq1}\frac{L\left(c+\epsilon\left(L\right)\right)}{\rho_{L,n}^2}\sum_{m\geq0}\left(\frac{L\left(c+\epsilon\left(L\right)\right)}{\rho_{L,n}}\right)^m-\frac{\eta}{\left(L+1\right)^2}\\
&=\frac{1}{L+1}\sum_{n\geq1}L\left(c+\epsilon\left(L\right)\right)\sum_{m\geq0}\frac{\left(L\left(c+\epsilon\left(L\right)\right)\right)^m}{\left(\rho_{L,n}\right)^{m+2}}-\frac{\eta}{\left(L+1\right)^2}\\
&=\frac{1}{L+1}\sum_{m\geq0}\left(L\left(c+\epsilon\left(L\right)\right)\right)^{m+1}\left(Z_\eta^{\left(2m+2\right)}+Z_\eta^{\left(2m+3\right)}\left(L\right)\right)-\frac{\eta}{\left(L+1\right)^2}
\end{align*}
or
\begin{equation}\label{main_eqn}
1=\frac{1}{L+1}\sum_{m\geq1}\left(L\left(c+\epsilon\left(L\right)\right)\right)^{m+1}\left(Z_\eta^{\left(2m\right)}\left(L\right)
+Z_\eta^{\left(2m+1\right)}\left(L\right)\right)-\frac{\eta\left(L\left(c+\epsilon\left(L\right)\right)\right)}{\left(L+1\right)^2},
\end{equation}
provided $L$ is sufficiently large. Now, let us write
\begin{equation}\label{epsilon_expansion}
	\epsilon\left(L\right)=\sum_{n=1}^{N-1}\frac{\epsilon_n}{L^n}+R_N\left(L\right),
\end{equation}
where the coefficient denoted as $\epsilon_n$ is determined by the recurrence relation \eqref{rec_reln}. Using the method of mathematical induction on the variable $N$, our aim to
prove that for all $N \geq 1$ and as $L \to \infty$, the relation $R_N(L) = \mathcal{O}_N(L^{-N})$ holds true. Throughout this paper, we utilize the subscript in the $\mathcal{O}$ notation to signify the dependence of the implied constant on certain parameters. It is important to note that this statement holds true for the case when $N=1$, since $R_1\left(L\right)=\epsilon\left(L\right)$. For $N \geq 2$ suppose that the statement holds for all $R_k\left(L\right)$ with $1\leq k \leq N-1$. Now, we are going to rewrite the right-hand side of equation \eqref{main_eqn}. First, we simplify
the right most term in equation \eqref{main_eqn} by using the induction hypothesis as follows
\begin{align*}
\frac{\eta L\left(c+\epsilon\left(L\right)\right)}{\left(L+1\right)^2}&=\frac{\eta}{L}\left(c+\epsilon\left(L\right)\right)\left(1+\frac{1}{L}\right)^{-2}\\
&=\frac{\eta}{L}\left(c+\sum_{n=1}^{N-2}\frac{\epsilon_n}{L^n}+R_{N-1}\left(L\right)\right)\sum_{n\geq0}\frac{\left(-1\right)^n\left(n+1\right)}{L^n}\\
&=\frac{\eta}{L}\left(c+\sum_{n=1}^{N-2}\frac{\epsilon_n}{L^n}+\mathcal{O}_{N-1}\left(\frac{1}{L^{N-1}}\right)\right)\sum_{n\geq0}\frac{\left(-1\right)^n\left(n+1\right)}{L^n}\\
&=\frac{c\eta}{L}\sum_{n\geq0}\left(-1\right)^n\frac{n+1}{L^n}+\frac{\eta}{L}\left(\sum_{n=1}^{N-2}\frac{\epsilon_n}{L^n}\right)\left(\sum_{n\geq0}\frac{\left(-1\right)^n\left(n+1\right)}{L^n}\right)+\mathcal{O}_N\left(\frac{1}{L^N}\right)\\
&=\frac{c \eta}{L}\sum_{n=0}^{N-2}\left(-1\right)^n\frac{n+1}{L^n}+\eta\sum_{n=0}^{N-3}\sum_{k=0}^{n}\left(-1\right)^{n-k}\left(n-k+1\right)\epsilon_{k+1}\frac{1}{L^{n+2}}+\mathcal{O}_N\left(\frac{1}{L^N}\right).
\end{align*}
Now we rewrite the term associated with zeta function $Z_\eta^{\left(2m\right)}\left(L\right)$ in equation \eqref{main_eqn} into three parts and analyze them separately:
\begin{align}\label{zeta_even}
\frac{1}{L+1}&\sum_{m\geq1}\left(L\left(c+\epsilon\left(L\right)\right)\right)^{m+1}Z_\eta^{\left(2m\right)}\left(L\right)\\
&=\frac{1}{L+1}\left(L\left(c+\epsilon\left(L\right)\right)\right)^2Z_\eta^{\left(2\right)}\left(L\right)+\frac{1}{L+1}
\sum_{m=2}^{N-1}\left(L\left(c+\epsilon\left(L\right)\right)\right)^{m+1}Z_\eta^{\left(2m\right)}\left(L\right)\nonumber\\
&+\frac{1}{L+1}\sum_{m\geq N}\left(L\left(c+\epsilon\left(L\right)\right)\right)^{m+1}Z_\eta^{\left(2m\right)}\left(L\right).\nonumber
\end{align}
Now, we analyze the first term of the right-hand side of equation \eqref{zeta_even} by using Lemma \ref{lemma_asy1} and we obtain
\begin{align*}
&\frac{1}{L+1}\left(L\left(c+\epsilon\left(L\right)\right)\right)^2Z_\eta^{\left(2\right)}\left(L\right)=
\frac{L}{L+1}\left(c+\epsilon\left(L\right)\right)^2\sum_{n\geq0}\frac{\zeta_{n,\eta}^{\left(2\right)}}{L^n}\\
&=\frac{1}{1+\frac{1}{L}}\left(c^2+\left(\epsilon\left(L\right)\right)^2+2c\epsilon\left(L\right)\right)\sum_{n\geq0}\frac{\zeta_{n,\eta}^{\left(2\right)}}{L^n}\\
&=\sum_{p\geq0}\frac{\left(-1\right)^p}{L^p}\sum_{n\geq0}\frac{\zeta_{n,\eta}^{\left(2\right)}}{L^n}\left(c^2+\left(\sum_{q=1}^{N-2}\frac{\epsilon_q}{L^q}
+\mathcal{O}_{N-1}\left(\frac{1}{L^{N-1}}\right)\right)^2+2c\left(\sum_{q=1}^{N-1}\frac{\epsilon_q}{L^q}+R_N\left(L\right)\right)\right)\\
&=\sum_{p\geq0}\frac{\left(-1\right)^p}{L^p}\sum_{n\geq0}\frac{\zeta_{n,\eta}^{\left(2\right)}}{L^n}\left(c^2+\left(\sum_{q=1}^{N-2}\frac{\epsilon_q}{L^q}\right)^2
+\mathcal{O}_{2N-2}\left(\frac{1}{L^{2N-2}}\right)+\mathcal{O}_{N}\left(\frac{1}{L^{N}}\right)+2c\left(\sum_{q=1}^{N-1}\frac{\epsilon_q}{L^q}+R_N\left(L\right)\right)\right)\\
&=\sum_{p\geq0}\frac{\left(-1\right)^p}{L^p}\sum_{n\geq0}\frac{\zeta_{n,\eta}^{\left(2\right)}}{L^n}\left(c^2+\left(\sum_{q=1}^{N-2}\frac{\epsilon_q}{L^q}\right)^2
+\mathcal{O}_{N}\left(\frac{1}{L^{N}}\right)+2c\sum_{q=1}^{N-1}\frac{\epsilon_q}{L^q}+2cR_N\left(L\right)\right)\\
&=\sum_{p\geq0}\frac{\left(-1\right)^p}{L^p}\sum_{n\geq0}\frac{\zeta_{n,\eta}^{\left(2\right)}}{L^n}\left(c^2+\left(\sum_{q=1}^{N-2}\frac{\epsilon_q}{L^q}\right)^2
+2c\sum_{q=1}^{N-1}\frac{\epsilon_q}{L^q}\right)+\mathcal{O}_{N}\left(\frac{1}{L^{N}}\right)+P\left(L\right)R_N\left(L\right)\\
&=\sum_{n=0}^{N-1}\sum_{k=0}^{n}\frac{\left(-1\right)^{n-k}\zeta_{k,\eta}^{\left(2\right)}}{L^n}\left(c^2+\left(\sum_{q=1}^{N-2}\frac{\epsilon_q}{L^q}\right)^2
+2c\sum_{q=1}^{N-1}\frac{\epsilon_q}{L^q}\right)+\mathcal{O}_{N}\left(\frac{1}{L^{N}}\right)+P\left(L\right)R_N\left(L\right),\\
\end{align*}
which can be written as
\begin{align*}
&c^2\sum_{n=0}^{N-1}\sum_{k=0}^{n}\frac{\left(-1\right)^{n-k}\zeta_{k,\eta}^{\left(2\right)}}{L^n}+\sum_{n=0}^{N-1}\sum_{k=0}^{n}\frac{\left(-1\right)^{n-k}\zeta_{k,\eta}^{\left(2\right)}}{L^n}
\sum_{l=0}^{N-3}\sum_{q=0}^{l}\frac{\epsilon_{q+1}\epsilon_{l-q+1}}{L^n}+2c\sum_{n=0}^{N-1}\sum_{k=0}^{n}\frac{\left(-1\right)^{n-k}\zeta_{k,\eta}^{\left(2\right)}}{L^n}\sum_{q=1}^{N-1}\frac{\epsilon_q}{L^q}\\
&=c^2\sum_{n=0}^{N-1}\sum_{k=0}^{n}\frac{\left(-1\right)^{n-k}\zeta_{k,\eta}^{\left(2\right)}}{L^n}+\sum_{n=0}^{N-3}\sum_{j=0}^{n}
\left(\sum_{k=0}^{n-j}\left(-1\right)^{n-j-k}\zeta_{k,\eta}^{\left(2\right)}\sum_{l=0}^{j}\epsilon_{l+1}\epsilon_{j-l+1}\right)\frac{1}{L^{n+2}}\\
&\quad\quad+2c\sum_{n=0}^{N-2}\sum_{k=0}^{n}\left(\epsilon_{k+1}\sum_{q=0}^{n-k}\left(-1\right)^{n-k-q}\zeta_{q,\eta}^{\left(2\right)}\right)\frac{1}{L^{n+1}}
+\mathcal{O}_N\left(\frac{1}{L^N}\right)+P\left(L\right)R_N\left(L\right),
\end{align*}
where $$P\left(L\right)=2c\sum_{p\geq0}\frac{\left(-1\right)^p}{L^p}\sum_{n\geq0}\frac{\zeta_{n,\eta}^{\left(2\right)}}{L^n}.$$
Next we analyze the second term on the right-hand side of equation \eqref{zeta_even} and obtain
\begin{align*}
\frac{1}{L+1}&\sum_{m=2}^{N-1}\left(L\left(c+\epsilon\left(L\right)\right)\right)^{m+1}Z_\eta^{\left(2m\right)}\left(L\right)\\
&=\left(1+\frac{1}{L}\right)^{-1}\sum_{m=2}^{N-1}\left(c+\sum_{n=1}^{N-m}\frac{\epsilon_n}{L^n}+R_{N-m+1}\left(L\right)\right)^{m+1}
\frac{1}{L^{m-1}}\sum_{k\geq0}\frac{\zeta_{k,\eta}^{\left(2m\right)}}{L^k}\\
&=\sum_{p\geq0}\frac{\left(-1\right)^p}{L^p}\sum_{m=2}^{N-1}\left(c+\sum_{n=1}^{N-m}\frac{\epsilon_n}{L^n}
+\mathcal{O}_{N-m+1}\left(\frac{1}{L^{N-m+1}}\right)\right)^{m+1}\frac{1}{L^{m-1}}\sum_{k\geq0}\frac{\zeta_{k,\eta}^{\left(2m\right)}}{L^k}\\
&=\sum_{p\geq0}\frac{\left(-1\right)^p}{L^p}\sum_{m=2}^{N-1}\left(\sum_{n=0}^{N-m}A_{m+1,n}\left(\epsilon_1,\ldots,\epsilon_n\right)\frac{1}{L^n}
+\mathcal{O}_{N-m+1}\left(\frac{1}{L^{N-m+1}}\right)\right)\frac{1}{L^{m-1}}\sum_{k\geq0}\frac{\zeta_{k,\eta}^{\left(2m\right)}}{L^k}
\end{align*}
or equivalently
\begin{align*}
&\sum_{p\geq0}\frac{\left(-1\right)^p}{L^p}\sum_{m=2}^{N-1}\frac{1}{L^{m-1}}\left(\sum_{n=0}^{N-m}\frac{1}{L^n}
\sum_{k=0}^{n}A_{m+1,k}\left(\epsilon_1,\ldots,\epsilon_k\right)\zeta_{n-k,\eta}^{\left(2m\right)}\right)+\mathcal{O}_{N}\left(\frac{1}{L^{N}}\right)\\
&=\sum_{n=0}^{N-3}\sum_{j=0}^{n}\left(-1\right)^{n-j}\sum_{m=2}^{j+2}\left(\sum_{k=0}^{j-m+2}\zeta_{j-m-k+2,\eta}^{2m}
A_{m+1,k}\left(\epsilon_1,\ldots,\epsilon_k\right)\right)\frac{1}{L^{n+1}}\\
&\quad\quad\quad+\sum_{j=0}^{N-2}\left(-1\right)^{N-j-2}\sum_{m=2}^{j+2}\left(\sum_{k=0}^{j-m+2}\zeta_{j-m-k+2,\eta}^{2m}
A_{m+1,k}\left(\epsilon_1,\ldots,\epsilon_k\right)\right)\frac{1}{L^{N-1}}+\mathcal{O}_N\left(\frac{1}{L^N}\right).
\end{align*}
Now, we analyze the last term on the right-hand side of \eqref{zeta_even} by using Lemma \ref{lemma_asy1}, Lemma \ref{Lemma3} and the induction hypothesis. In this case we have that
\begin{align*}
&\frac{1}{L+1}\sum_{m\geq N}\left(L\left(c+\epsilon\left(L\right)\right)\right)^{m+1}Z_\eta^{\left(2m\right)}\left(L\right)\\
&=\frac{1}{L+1}\left(\left(L\left(c+\epsilon\left(L\right)\right)\right)^{N+1}Z_\eta^{\left(2N\right)}\left(L\right)+\sum_{m\geq N+1}  \left(L\left(c+\epsilon\left(L\right)\right)\right)^{m+1}Z_\eta^{\left(2m\right)}\left(L\right) \right)\\
&=\frac{1}{1+\frac{1}{L}}L^N\left(c+\epsilon\left(L\right)\right)^{N+1}Z_\eta^{\left(2N\right)}\left(L\right)+\frac{1}{L+1}\sum_{m\geq N+1}  \left(L\left(c+\epsilon\left(L\right)\right)\right)^{m+1}Z_\eta^{\left(2m\right)}\left(L\right),
\end{align*}
which can be written as
\begin{align*}
&\frac{1}{1+\frac{1}{L}}\frac{\left(c+\epsilon\left(L\right)\right)^{N+1}}{L^{N-1}}\sum_{n\geq0}\frac{\zeta_{n,\eta}^{\left(2N\right)}}{L^n}+\frac{1}{L+1}\sum_{m\geq N+1}  \left(L\left(c+\epsilon\left(L\right)\right)\right)^{m+1}Z_\eta^{\left(2m\right)}\left(L\right)\\
&=\frac{1}{1+\frac{1}{L}}\frac{\left(c+\epsilon\left(L\right)\right)^{N+1}}{L^{N-1}}\left(\zeta_{0,\eta}^{\left(2N\right)}+\sum_{n\geq1}\frac{\zeta_{n,\eta}^{\left(2N\right)}}{L^n} \right)+\frac{1}{L+1}\sum_{m\geq N+1}\left(L\left(c+\epsilon\left(L\right)\right)\right)^{m+1}Z_\eta^{\left(2m\right)}\left(L\right)\\
&=\frac{c^{N+1}\zeta_{0,\eta}^{2N}}{L^{N-1}}+\mathcal{O}_N\left(\frac{1}{L^N}\right)+\frac{1}{L+1}\sum_{m\geq N+1}\left(L\left(c+\epsilon\left(L\right)\right)\right)^{m+1}Z_\eta^{\left(2m\right)}\left(L\right)\\
&=\frac{c^{N+1}\zeta_{0,\eta}^{2N}}{L^{N-1}}+\mathcal{O}_N\left(\frac{1}{L^N}\right)+\frac{1}{1+\frac{1}{L}}\sum_{m\geq N+1}L^m \left(c+\epsilon\left(L\right)\right)^{m+1}Z_\eta^{\left(2m\right)}\left(L\right)\\
&=\frac{c^{N+1}\zeta_{0,\eta}^{2N}}{L^{N-1}}+\mathcal{O}_N\left(\frac{1}{L^N}\right)+\frac{1}{1+\frac{1}{L}}\mathcal{O}\left(1\right)\sum_{m\geq N+1}\frac{\left(c+\epsilon\left(L\right)\right)^{m+1}}{L^{m-1}}\\
&=\frac{c^{N+1}\zeta_{0,\eta}^{2N}}{L^{N-1}}+\mathcal{O}_N\left(\frac{1}{L^N}\right)+\frac{1}{1+\frac{1}{L}}\mathcal{O}\left(1\right)\frac{\left(c+\epsilon\left(L\right)\right)^{N+2}}{L^N}
\sum_{m\geq0}\frac{\left(c+\epsilon\left(L\right)\right)^{m}}{L^{m}}\\
&=\frac{c^{N+1}\zeta_{0,\eta}^{2N}}{L^{N-1}}+\mathcal{O}_N\left(\frac{1}{L^N}\right)+\frac{1}{1+\frac{1}{L}}\mathcal{O}\left(1\right)\frac{\left(c+o\left(1\right)\right)^{N+2}}{L^N}
\sum_{m\geq0}\frac{\left(c+o\left(1\right)\right)^{m}}{L^{m}}\\
&=\frac{c^{N+1}\zeta_{0,\eta}^{2N}}{L^{N-1}}+\mathcal{O}_N\left(\frac{1}{L^N}\right)+\mathcal{O}\left(\frac{1}{L^N}\right)=\frac{c^{N+1}\zeta_{0,\eta}^{2N}}{L^{N-1}}+\mathcal{O}_N\left(\frac{1}{L^N}\right).
\end{align*}
Now we divide the term associated with $Z_\eta^{\left(2m+1\right)}\left(L\right)$ in equation \eqref{main_eqn} into two parts and analyze them separately. To do this observe that the next relation holds true:
\begin{align}\label{zeta_odd}
&\frac{1}{L+1}\sum_{m\geq1}\left(L\left(c+\epsilon\left(L\right)\right)\right)^{m+1}Z_\eta^{\left(2m+1\right)}\left(L\right)\\
&=\frac{1}{L+1}\sum_{m=1}^{N-3}\left(L\left(c+\epsilon\left(L\right)\right)\right)^{m+1}Z_\eta^{\left(2m+1\right)}\left(L\right)
+\frac{1}{L+1}\sum_{m\geq N-2}\left(L\left(c+\epsilon\left(L\right)\right)\right)^{m+1}Z_\eta^{\left(2m+1\right)}\left(L\right).\nonumber
\end{align}
First, we analyze the first term on the right-hand side of equation \eqref{zeta_odd} by using the induction hypothesis and Lemma \ref{lemma_asy1}. Observe that this term
\begin{align*}
\frac{1}{L+1}&\sum_{m=1}^{N-3}\left(L\left(c+\epsilon\left(L\right)\right)\right)^{m+1}Z_\eta^{\left(2m+1\right)}\left(L\right)\\
&=\left(1+\frac{1}{L}\right)^{-1}\sum_{m=1}^{N-3}\left(c+\epsilon\left(L\right)\right)^{m+1}\frac{1}{L^{m+1}}\sum_{k\geq0}\frac{\zeta_{k,\eta}^{\left(2m+1\right)}}{L^k}\\
&=\sum_{p\geq0}\frac{\left(-1\right)^p}{L^p}\sum_{m=1}^{N-3}\left(c+\sum_{n=1}^{N-m-2}\frac{\epsilon_n}{L^n}+R_{N-m-1}\left(L\right)\right)^{m+1}\frac{1}{L^{m+1}}
\sum_{k\geq0}\frac{\zeta_{k,\eta}^{\left(2m+1\right)}}{L^k}\\
&=\sum_{p\geq0}\frac{\left(-1\right)^p}{L^p}\sum_{m=1}^{N-3}\left(c+\sum_{n=1}^{N-m-2}\frac{\epsilon_n}{L^n}
+\mathcal{O}_{N-m-1}\left(\frac{1}{L^{N-m-1}}\right)\right)^{m+1}\frac{1}{L^{m+1}}\sum_{k\geq0}\frac{\zeta_{k,\eta}^{\left(2m+1\right)}}{L^k}\\
&=\sum_{p\geq0}\frac{\left(-1\right)^p}{L^p}\sum_{m=1}^{N-3}\left(\sum_{n=0}^{N-m-2}A_{m+1,n}\left(\epsilon_1,\ldots,\epsilon_n\right)\frac{1}{L^n}
+\mathcal{O}_{N-m-1}\left(\frac{1}{L^{N-m-1}}\right)\right)\frac{1}{L^{m+1}}\sum_{k\geq0}\frac{\zeta_{k,\eta}^{\left(2m+1\right)}}{L^k}
\end{align*}
can be written as
\begin{align*}
&\sum_{p\geq0}\frac{\left(-1\right)^p}{L^p}\sum_{m=1}^{N-3}\frac{1}{L^{m+1}}\left(\sum_{n=0}^{N-m-2}A_{m+1,n}\left(\epsilon_1,\ldots,\epsilon_n\right)\frac{1}{L^n}
\sum_{k\geq0}\frac{\zeta_{k,\eta}^{\left(2m+1\right)}}{L^k}\right)+\mathcal{O}_{N}\left(\frac{1}{L^{N}}\right)\\
&=\sum_{p\geq0}\frac{\left(-1\right)^p}{L^p}\sum_{m=1}^{N-3}\frac{1}{L^{m+1}}\left(\sum_{n=0}^{N-m-2}\frac{1}{L^n}
\sum_{k=0}^{n}A_{m+1,k}\left(\epsilon_1,\ldots,\epsilon_k\right)\zeta_{n-k,\eta}^{\left(2m+1\right)}\right)+\mathcal{O}_{N}\left(\frac{1}{L^{N}}\right)\\
&=\sum_{n=0}^{N-4}\sum_{j=0}^{n}\sum_{m=1}^{j+1}\sum_{k=0}^{j-m+1}\left(-1\right)^{n-j}\zeta_{j-m-k+1,\eta}^{2m+1}
A_{m+1,k}\left(\epsilon_1,\ldots,\epsilon_k\right)\frac{1}{L^{n+2}}\\
&\quad\quad\quad+\sum_{j=0}^{N-3}\sum_{m=1}^{j+1}\sum_{k=0}^{j-m+1}\left(-1\right)^{N-j-3}\zeta_{j-m-k+1,\eta}^{2m+1}A_{m+1,k}\left(\epsilon_1,\ldots,\epsilon_k\right)\frac{1}{L^{N-1}}+\mathcal{O}_N\left(\frac{1}{L^N}\right).
\end{align*}
Now, we focus on the second term on the right-hand side of equation \eqref{zeta_odd} by employing the Lemma \ref{lemma_asy1}, Lemma \ref{Lemma3} as well as the induction hypothesis and obtain
\begin{align*}
&\frac{1}{L+1}\sum_{m\geq N-2}\left(L\left(c+\epsilon\left(L\right)\right)\right)^{m+1}Z_\eta^{\left(2m+1\right)}\left(L\right)\\
&=\frac{1}{L+1}\left(\left(L\left(c+\epsilon\left(L\right)\right)\right)^{N-1}Z_\eta^{\left(2N-3\right)}\left(L\right)+\sum_{m\geq N-1}  \left(L\left(c+\epsilon\left(L\right)\right)\right)^{m+1}Z_\eta^{\left(2m+1\right)}\left(L\right) \right)\\
&=\frac{1}{1+\frac{1}{L}}L^{N-2}\left(c+\epsilon\left(L\right)\right)^{N-1}Z_\eta^{\left(2N-3\right)}\left(L\right)+\frac{1}{L+1}\sum_{m\geq N-1}  \left(L\left(c+\epsilon\left(L\right)\right)\right)^{m+1}Z_\eta^{\left(2m+1\right)}\left(L\right)\\
&=\frac{1}{1+\frac{1}{L}}\frac{\left(c+\epsilon\left(L\right)\right)^{N-1}}{L^{N-1}}\sum_{n\geq0}\frac{\zeta_{n,\eta}^{\left(2N-3\right)}}{L^n}+\frac{1}{L+1}\sum_{m\geq N-1}  \left(L\left(c+\epsilon\left(L\right)\right)\right)^{m+1}Z_\eta^{\left(2m+1\right)}\left(L\right)\\
&=\frac{1}{1+\frac{1}{L}}\frac{\left(c+\epsilon\left(L\right)\right)^{N-1}}{L^{N-1}}\left(\zeta_{0,\eta}^{\left(2N-3\right)}+\sum_{n\geq1}\frac{\zeta_{n,\eta}^{\left(2N-3\right)}}{L^n} \right)
+\frac{1}{L+1}\sum_{m\geq N-1}\left(L\left(c+\epsilon\left(L\right)\right)\right)^{m+1}Z_\eta^{\left(2m+1\right)}\left(L\right)\\
&=\frac{c^{N-1}\zeta_{0,\eta}^{2N-3}}{L^{N-1}}+\mathcal{O}_N\left(\frac{1}{L^N}\right)+\frac{1}{L+1}\sum_{m\geq N-1}\left(L\left(c+\epsilon\left(L\right)\right)\right)^{m+1}Z_\eta^{\left(2m+1\right)}\left(L\right)
\end{align*}
or equivalently
\begin{align*}
&\frac{c^{N-1}\zeta_{0,\eta}^{2N-3}}{L^{N-1}}+\mathcal{O}_N\left(\frac{1}{L^N}\right)+\frac{1}{\left(1+\frac{1}{L}\right)}\sum_{m\geq N-1}L^m \left(c+\epsilon\left(L\right)\right)^{m+1}Z_\eta^{\left(2m+1\right)}\left(L\right)\\
&=\frac{c^{N-1}\zeta_{0,\eta}^{2N-3}}{L^{N-1}}+\mathcal{O}_N\left(\frac{1}{L^N}\right)+\frac{1}{\left(1+\frac{1}{L}\right)}\mathcal{O}\left(1\right)\sum_{m\geq N-1}  \frac{\left(c+\epsilon\left(L\right)\right)^{m+1}}{L^{m+1}}\\
&=\frac{c^{N-1}\zeta_{0,\eta}^{2N-3}}{L^{N-1}}+\mathcal{O}_N\left(\frac{1}{L^N}\right)+\frac{1}{\left(1+\frac{1}{L}\right)}\mathcal{O}\left(1\right)\frac{\left(c+\epsilon\left(L\right)\right)^{N}}{L^N}\sum_{m\geq0}  \frac{\left(c+\epsilon\left(L\right)\right)^{m}}{L^{m}}\\
&=\frac{c^{N-1}\zeta_{0,\eta}^{2N-3}}{L^{N-1}}+\mathcal{O}_N\left(\frac{1}{L^N}\right)+\frac{1}{\left(1+\frac{1}{L}\right)}\mathcal{O}\left(1\right)  \frac{\left(c+o\left(1\right)\right)^{N}}{L^N}\sum_{m\geq0}\frac{\left(c+o\left(1\right)\right)^{m}}{L^{m}}\\
&=\frac{c^{N-1}\zeta_{0,\eta}^{2N-3}}{L^{N-1}}+\mathcal{O}_N\left(\frac{1}{L^N}\right)+\mathcal{O}\left(\frac{1}{L^N}\right)=\frac{c^{N-1}\zeta_{0,\eta}^{2N-3}}{L^{N-1}}+\mathcal{O}_N\left(\frac{1}{L^N}\right).
\end{align*}
Substituting these results into \eqref{main_eqn} and equating the constant term and the coefficient of $\frac{1}{L}$, we obtain that
$c=\sqrt{2}$ and $$2c\zeta_{0,\eta}^{\left(2\right)}\epsilon_1=c\eta-c^2\sum_{k=0}^{1}\left(-1\right)^{1-k}\zeta_{k,\eta}^{\left(2\right)}-\zeta_{0,\eta}^{\left(4\right)}A_{3,0}.$$ Moreover, applying the recurrence relation \eqref{rec_reln} in the form
\begin{align*}
&\left(-1\right)^{n}c\eta \left(n+2\right)+\eta\sum_{k=0}^{n}\left(-1\right)^{n-k+1}\left(n-k+1\right)\epsilon_{k+1}+c^2\sum_{k=0}^{n+2}\left(-1\right)^{n-k+2}\zeta_{k,\eta}^{\left(2\right)}\\
&\quad\quad+\sum_{j=0}^{n}\left(\sum_{k=0}^{n-j}\left(-1\right)^{n-j-k}\zeta_{k,\eta}^{\left(2\right)}\sum_{l=0}^{j}\epsilon_{l+1}\epsilon_{j-l+1}\right)\\
&\quad\quad+2c\sum_{k=0}^{n+1}\epsilon_{k+1}\sum_{q=0}^{n-k+1}\left(-1\right)^{n-k-q+1}\zeta_{q,\eta}^{\left(2\right)}+\sum_{j=0}^{n+1}\left(-1\right)^{n-j+1}
\sum_{m=2}^{j+2}\left(\sum_{k=0}^{j-m+2}\zeta_{j-m-k+2,\eta}^{\left(2m\right)}A_{m+1,k}\left(\epsilon_1,\ldots,\epsilon_k\right)\right)\\
&\quad\quad+\sum_{j=0}^{n}\left(-1\right)^{n-j}\sum_{m=1}^{j+1}\left(\sum_{k=0}^{j-m+1}\zeta_{j-m-k+1,\eta}^{\left(2m+1\right)}A_{m+1,k}\left(\epsilon_1,\ldots,\epsilon_k\right)\right)=0,
\end{align*}
equation \eqref{main_eqn} simplifies to
$$1=\mathcal{O}_N\left(\frac{1}{L^N}\right)+\mathcal{O}_{N}\left(\frac{1}{L^{N}}\right)+P\left(L\right)R_N\left(L\right)+\mathcal{O}_N\left(\frac{1}{L^N}\right)
+\mathcal{O}_N\left(\frac{1}{L^N}\right)+\mathcal{O}_N\left(\frac{1}{L^N}\right)+\mathcal{O}_N\left(\frac{1}{L^N}\right),$$
that is,
$$R_N\left(L\right)P\left(L\right)=\mathcal{O}_N\left(\frac{1}{L^N}\right)$$
as $L\to +\infty$. Since $$P\left(L\right)=2c\sum_{p\geq0}\frac{\left(-1\right)^p}{L^p}\sum_{n\geq0}\frac{\zeta_{n,\eta}^{\left(2\right)}}{L^n},$$ it follows that $R_N\left(L\right)=\mathcal{O}_N\left(\frac{1}{L^N}\right)$, which completes the proof of Theorem \ref{asym_the}.
\end{proof}

\subsection*{Acknowledgements} Pranav Kumar is grateful to the Council of Scientific and Industrial Research India (Grant No. 09/1022(0060)/2018-EMR-I) for the financial support, and Sanjeev Singh is thankful to the Science and Engineering Research Board (SERB), Department of Science and Technology, Government of India for the financial support through Project MTR/2022/000792.

\end{document}